\documentclass[11pt]{amsart} 
\usepackage{amssymb,psfig}

\DeclareMathOperator{\dom}{dom} 
\newcommand{\dbar}{\ensuremath{\overline\partial}} 
\newcommand{\dbarstar}{\ensuremath{\overline\partial^*}} 
\newcommand{\C}{\ensuremath{\mathbb{C}}} 
\newcommand{\R}{\ensuremath{\mathbb{R}}} 
 
\newcommand{\T}{\ensuremath{\mathbb{T}}}

\begin{document} 
 
\numberwithin{equation}{section} 
 
\newtheorem{theorem}{Theorem}[section] 
\newtheorem{proposition}[theorem]{Proposition} 
\newtheorem{conjecture}[theorem]{Conjecture} 
\def\theconjecture{\unskip} 
\newtheorem{corollary}[theorem]{Corollary} 
\newtheorem{lemma}[theorem]{Lemma} 
\newtheorem{observation}[theorem]{Observation} 
\theoremstyle{definition} 
\newtheorem{definition}{Definition} 
\numberwithin{definition}{section} 
\newtheorem{remark}{Remark} 
\def\theremark{\unskip} 
\newtheorem{question}{Question} 
\def\thequestion{\unskip} 
\newtheorem{example}{Example} 
\def\theexample{\unskip} 
\newtheorem{problem}{Problem} 
 
\def\intprod{\mathbin{\lr54}} 
\def\reals{{\mathbb R}} 
\def\integers{{\mathbb Z}} 
\def\naturals{{\mathbb N}} 
\def\complex{{\mathbb C}\/} 
\def\distance{\operatorname{distance}\,} 
\def\sgn{\operatorname{sgn}\,} 
\def\ZZ{ {\mathbb Z} } 
\def\e{\varepsilon} 
\def\p{\partial} 
\def\rp{{ ^{-1} }} 
\def\Re{\operatorname{Re\,} } 
\def\Im{\operatorname{Im\,} } 
\def\dbarb{\bar\partial_b} 
\def\eps{\varepsilon}

\def\scriptu{{\mathcal U}} 
\def\scriptr{{\mathcal R}} 
\def\scripta{{\mathcal A}} 
\def\scripti{{\mathcal I}} 
\def\scripth{{\mathcal H}} 
\def\scriptm{{\mathcal M}} 
\def\scripte{{\mathcal E}} 
\def\scriptt{{\mathcal T}} 
\def\scriptb{{\mathcal B}} 
\def\scriptf{{\mathcal F}} 
\def\scripto{{\mathfrak o}} 
\def\scriptv{{\mathcal V}} 
\def\frakg{{\mathfrak g}} 
\def\frakG{{\mathfrak G}} 
 
\def\ov{\overline} 
\date {April 21, 2003. Minor revision November 5, 2003.}
 
\author{Michael Christ}
\thanks{M.C.\ was supported in part by NSF\ grant DMS 9970660.}
\address{University of California, Berkeley}
\email{mchrist@math.berkeley.edu}

\author{Siqi Fu} 
\thanks
{S.F.\ was supported in part by NSF
grant DMS 0070697 and by an AMS centennial fellowship.}
\address{University of Wyoming}
\email{sfu@uwyo.edu}
\title[Compactness in the $\dbar$-Neumann problem and the Aharonov-Bohm effect]  
{ Compactness in the $\dbar$-Neumann problem, \\Magnetic Schr\"odinger operators, 
\\ and the Aharonov-Bohm effect} 

\begin{abstract}
Compactness of the Neumann operator in the d-bar Neumann problem is studied
for weakly pseudoconvex bounded Hartogs domains in two dimensions. A nonsmooth
example is given in which condition (P) fails to hold, yet the Neumann operator
is compact. The main result, in contrast, is that for smoothly bounded Hartogs domains,
condition (P) of Catlin and Sibony is equivalent to compactness. 

The analyses of both compactness and condition (P) boil down to properties of the lowest
eigenvalues of certain sequences of Schrodinger operators, with and without magnetic fields, 
parametrized by a Fourier variable resulting from the Hartogs symmetry. The nonsmooth
counterexample is based on the Aharonov-Bohm phenomenon of quantum mechanics. 
For smooth domains, we prove that there always exists an exceptional sequence of Fourier variables
for which the Aharonov-Bohm effect is weak. This sequence can be very sparse, so
that the failure of compactness is due to a rather subtle effect.
\end{abstract}

\maketitle

 
 
\section{Introduction} 
Let $\Omega$ be a bounded domain in $\complex^n$.  The 
$\dbar$-Neumann Laplacian $\square=\dbar\dbarstar +\dbarstar\dbar$ 
is a formally self-adjoint operator acting on $(0, q)$-forms with 
$L^2$-coefficients satisfying certain boundary conditions.  
The Kohn Laplacian $\square_b$ is a 
non-elliptic operator acting on forms on the boundary, defined 
under certain regularity assumptions on $b\Omega$~\cite{Kohn65}.  
An extensive literature is devoted to the problem of 
relating complex-geometric properties of $\partial\Omega$ 
with analytical properties of the $\dbar$-Neumann problem and $\square_b$. 
Kohn ~\cite{Kohn6364} analyzed the $\dbar$-Neumann 
problem for smoothly bounded strictly pseudoconvex domains, 
and subsequently the subelliptic theory of  
$\square$ and $\square_b$ has become in large part 
understood \cite{Catlin83, Catlin87, Dangelo82, Kohn79}. 
However, various fundamental issues for domains of infinite type, 
for which no subelliptic estimates hold, remain unresolved. 
See for instance \cite{boasstraube99},\cite{Christ99},\cite{Dangelokohn99}
for surveys of aspects of the $\bar\partial$-Neumann problem and Kohn Laplacian.
Some recent work on global regularity and on $C^\infty$ hypoellipticity, 
for domains of infinite type, is in \cite{Christ96, Christ01a, Christ01b}.
 
In this paper, we study compactness of $\square_b$ and the $\dbar$-Neumann problem. 
For smoothly bounded pseudoconvex domains, it is well-known that 
compactness is a property weaker than subellipticity, but 
stronger than $C^\infty$ global regularity~\cite{KohnNirenberg65}.  
The well-known property ($P$) was first 
introduced by Catlin, who proved that it 
implies compactness for smoothly bounded pseudoconvex 
domains~\cite{Catlin84}, and showed that it is implied by natural geometric
conditions. 
It was later systematically 
studied by Sibony for all compact sets in $\C^n$ from the viewpoint 
of potential theory \cite{Sibony87}.  A compact set 
$K$ in $\complex^n$ is said to satisfy property ($P$) (or to be $B$-regular in 
the terminology of Sibony) if for any $M>0$, there exist a 
neighborhood $U$ of $K$ and a function $\rho\in C^\infty(U)$, 
$0\le\rho\le 1$, such that the complex Hessian $(\partial^2 \rho 
/\partial z_j\partial\bar z_k)$ is $\ge M$ at every point of $U$.   
Straube proved that 
Catlin's result on the $\dbar$-Neumann Laplacian holds for all 
bounded pseudoconvex domains without any regularity assumption on 
the boundary~\cite{Straube97}.  
 
It has long been known that 
compactness precludes the presence of complex discs in the boundary for domains 
in $\C^2$ (under minimal regularity assumptions on the boundary, 
say Lipschitz). The converse is not true; 
Matheos~\cite{Matheos97} constructed a smoothly 
bounded, pseudoconvex, complete Hartogs domain in $\C^2$ whose 
boundary contains no complex analytic disc but whose 
$\dbar$-Neumann Laplacian nevertheless does not have compact 
resolvent. (See~\cite{FuStraube01} for a discussion of this and 
other results on compactness.)   However, whether compactness is 
equivalent to property ($P$) for smoothly bounded pseudoconvex 
domains in $\C^2$ had remained an open question\footnote{
McNeal \cite{Mcneal02} has introduced a variant ($\tilde P$)
of condition ($P$). ($P$) implies ($\tilde P$) for all smooth pseudoconvex
domains, and the two are equivalent for Hartogs domains, which are
the only domains discussed in this paper.}, which we explore and
answer for Hartogs domains --- both in the affirmative and in the negative --- in this paper.   
 
Matheos \cite{Matheos97} exploited the equivalence between
compactness in the $\dbar$-Neumann problem, 
for smoothly bounded Hartogs domains in $\complex^2$, 
and a certain property of an associated
one-parameter family of magnetic Schr\"odinger operators in $\complex^1$.
Fu and Straube \cite{FuStraube02} observed that for smoothly bounded, pseudoconvex, 
complete Hartogs domains in $\C^2$,  this problem is closely related to 
topics discussed in the mathematical physics literature under the names
{\em diamagnetism} and {\em paramagnetism}. 
More precisely, that property ($P$) implies compactness for Hartogs
domains is a consequence of diamagnetism, and whether compactness implies property ($P$) 
is connected to paramagnetism.
For more on diamagnetism and paramagnetism the reader
may consult \cite{Simon76},\cite{Simon79},\cite{AvronSimon79},\cite{Erdos97}.
 
Our results for compactness of the $\dbar$-Neumann problem and
Kohn Laplacian are twofold. By a {\em Hartogs domain} we mean
an open subset $\Omega\subset \complex^{d+1}$ such that
whenever $(z,w)\in\Omega$, likewise $(z,e^{i\theta}w)\in\Omega$
for every $\theta\in\reals$.
Such a domain is said to be complete if whenever $(z,w)\in\Omega$
and $|w'|\le|w|$, $(z,w')\in\Omega$.

\begin{theorem} \label{thm:dbarsmoothcase}
Let $\Omega\subset\C^2$ be a
smoothly bounded pseudoconvex Hartogs domain. 
Then the following are equivalent:

\begin{enumerate}
\item The $\dbar$-Neumann Laplacian $\square$ has compact
resolvent in $L^2(\Omega)$.
\item The Kohn Laplacian $\square_b$ has compact resolvent.
\item $b\Omega$ satisfies property ($P$).
\end{enumerate}
\end{theorem}
\noindent For the precise meaning of (2) see Definition~\ref{compactnessdefn} below.
 
The equivalence between compactness and property ($P$) is however a quantitative 
rather than a qualitative phenomenon, which breaks down for boundaries having
very limited regularity.
\begin{theorem}\label{compactness} 
There exists a pseudoconvex, complete Hartogs domain 
$\Omega=\{(z, w)\in \complex^2: \ |w|< e^{-\psi(z)}, \ |z|<2\}$, where $\psi$ 
is continuous, $\nabla\psi\in L^2$ in any compact subset of $\{|z|<2\}$,  
and $\Delta\psi\in L^1$ is lower semicontinuous, 
such that $b\Omega$ does not satisfy property {\rm (P)}, yet  
a Kohn Laplacian is well-defined on $b\Omega$ and satisfies a compactness inequality.
\end{theorem} 
However, we must emphasize the distinction between 
``a Kohn Laplacian'' and ``the Kohn Laplacian''. The latter 
is defined with respect to the Hilbert space structure 
$L^2(b\Omega)$ induced by surface measure on $b\Omega$, 
while the Kohn Laplacian of our theorem is defined in terms  
of a different measure, which is smooth when  
expressed in terms of certain natural coordinates for $b\Omega$, 
but is quite different from surface measure.  
Our example could alternatively be described as a nonsmooth 
three dimensional CR manifold with an $S^1$ action, 
or as the unit sphere bundle in a holomorphic line bundle over a 
one-dimensional base manifold, equipped with a nonsmooth 
metric. 
 
It is well-known that Schr\"odinger operators with magnetic fields arise 
in connection with holomorphic line bundles over complex manifolds; see
for instance \cite{demailly}.
They arise in connection with Hartogs domains for the same reasons.
Our results both amount to a semi-classical analysis of certain 
magnetic Schr\"{o}dinger operators.  Let $\varphi$ be a 
subharmonic function on the unit disc $\Omega_0$ such that 
$\nabla\varphi\in L^2(\Omega_0)$ in the sense of distributions. 
Let $S_\varphi$ be a Schr\"{o}dinger operator formally given by\footnote{
This is essentially a Pauli operator. The two-dimensional Pauli operator
$\sigma\cdot(\nabla-a)^2$, with $a = i(\varphi_y,-\varphi_x)$,
splits into two direct summands, $S_\varphi-\Delta\varphi$ and
$S_{\varphi}+\Delta\varphi$. Nonnegativity of $\Delta\varphi$
implies that $S_\varphi-\Delta\varphi\le S_\varphi
\le 2(S_\varphi-\Delta\varphi)$.}
\begin{equation*} 
S_\varphi = -[(\p_x+i\varphi_y)^2 + (\p_y -i\varphi_x)^2] + 
\Delta\varphi 
\end{equation*} 
and let $S^0_{\varphi} = -\p_x^2-\p_y^2 + \Delta\varphi$ be the 
Schr\"{o}dinger operator with the same electric potential but zero 
magnetic potential.  Denote by $\lambda^{\rm m}_\varphi$ and 
$\lambda^{\rm e}_\varphi$ respectively the lowest eigenvalues of the Dirichlet 
realizations of $S_\varphi$ and of $S^0_\varphi$.
The diamagnetic inequality of Simon \cite{Simon76} (see also Kato~\cite{Kato72} and
Simon~\cite{Simon79}) guarantees that $\lambda^{\rm e}_\varphi
\le \lambda^{\rm m}_\varphi$ for any $\varphi$.
For $C^\infty$ pseudoconvex complete Hartogs domains which are strictly pseudoconvex
in a neighborhood of $\{(z,w): w=0\}$,
Fu and Straube \cite{FuStraube02} proved that condition ($P$) is equivalent to
$\lambda_{n\varphi}^{\rm e}\to \infty$ as $n\to+\infty$. It is implicit in
the analysis of Matheos \cite{Matheos97} that $\lambda_{n\varphi}^{\rm m}\to \infty$
is equivalent to compactness.

It follows from \cite{Sibony87} 
that there exists a $C^\infty$, pseudoconvex, complete Hartogs domain
$\Omega\subset\complex^2$ for which the set of weakly pseudoconvex boundary points
has positive measure, yet $b\Omega$ satisfies property ($P$); and there exists
another such domain for which the set of strictly pseudoconvex points is dense,
yet $b\Omega$ does not satisfy property ($P$). In light of this and  \cite{FuStraube02},
there exists a $C^\infty$ subharmonic function $\psi$ for which $\Delta\psi=0$
on a set of positive Lebesgue measure, yet
$\lambda^e_{n\psi}\to\infty$; and on the other hand there exists another such $\psi$
for which $\{z:\Delta\psi>0\}$ is dense, yet $\sup_n \lambda^e_{n\psi}<\infty$.

By virtue of these equivalences, Theorem~\ref{compactness} amounts to:
\begin{theorem}\label{main} 
There exists a continuous subharmonic function $\varphi$ on the unit disk 
$\Omega_0$ with $\nabla\varphi\in L^2 (\Omega_0)$ in the sense of 
distributions and $\Delta\varphi\in L^1(\Omega_0)$ lower semicontinuous, such 
that $\lim_{n\to\infty}\lambda^{\rm m}_{n\varphi}=\infty$ but 
$\lim_{n\to\infty}\lambda^{\rm e}_{n\varphi}<\infty$. 
\end{theorem} 

This degree of regularity is quite natural from the perspective of Schr\"odinger 
operators, as it guarantees that the magnetic Schr\"odinger form 
$Q(u,v)=\langle S_\varphi u,\,v\rangle$ is well-defined for all $u,v$ in a standard dense 
subclass of $L^2(\Omega_0)$, namely $C^1_0(\Omega_0)$. 
 
Theorem~\ref{main} is based on the
Aharonov-Bohm effect \cite{aharonovbohm}, a quantum phenomenon
in which a physical system not exposed to
a magnetic {\em field} is nonetheless influenced by the associated magnetic {\em potential}.
Avron and Simon~\cite{AvronSimon79} gave a counterexample, based in part on this effect,
to a conjectured paramagnetic inequality, which when specialized to our situation would
have implied that $\lambda^{\rm m}_{n\varphi}$ is always $\le \lambda^{\rm e}_{2n\varphi}$. 
Theorem~\ref{main} realizes a more extreme form of this phenomenon,
providing an example where paramagnetism can fail more dramatically.
 
The following weaker variant for $C^\infty$ structures is 
an easy consequence of a simpler form of the same construction. 
 
\begin{proposition} \label{prop:smoothcase}
There exists a $C^\infty$ subharmonic function $\varphi$ on the unit disk 
$\Omega_0$  
such that 
$\sup_{n\to\infty}\lambda^{\rm e}_{n\varphi}<\infty$ 
but
$\limsup_{n\to\infty}\lambda^{\rm m}_{n\varphi}=\infty$.
\end{proposition} 
 
One can even make $\lambda^{\rm m}_{n \varphi}\to\infty$ as $n\to\infty$ 
through a subset of $\naturals$ whose complement is quite sparse. 
But to control every value of $n$, without exception, is a different matter. 

\begin{theorem}  \label{thm:smoothcase}
Let $\varphi$ be subharmonic, and suppose that 
$\Delta\varphi$ is H\"older continuous of some positive order.
If $\sup_n\lambda^{\rm e}_{n\varphi}<\infty$ then
$\liminf_{n\to\infty}\lambda^{\rm m}_{n\varphi}<\infty$.
\end{theorem}

Our analysis produces a concrete bound for the rate of growth of a subsequence $(n_\nu)$ 
for which $\lambda^{\rm m}_{n_\nu\varphi}$ remains bounded,
but this bound allows for sequences of very large gaps and strongly suggests that
there should exist domains for which such subsequence all have upper density zero. 
The pigeonhole principle plays a crucial role in the proof that such a subsequence must
exist.
Theorem~\ref{thm:dbarsmoothcase} is a consequence of Theorem~\ref{thm:smoothcase};
it is only through the existence of possibly sparse sequences of exceptional
values of $n$ that the failure of property ($P$) implies the failure of compactness.

Although our work may be viewed as a semiclassical analysis of magnetic
Schr\"odinger operators, the point of view is different than that ordinarily
taken in mathematical physics. There one studies $(h\nabla-A)^2$ as
$h$ tends to zero. We instead analyze $(\nabla -h\rp A)^2$, as $h\to 0$,  and
are interested in whether the lowest eigenvalue tends to infinity. 
In semiclassical terms we have a situation where the lowest eigenvalue is 
positive and tends to zero with $h$, and we are interested in whether or not
it is $O(h^2)$. The other distinction is that our magnetic field is an arbitrary
nonnegative function (with certain regularity), rather than a function with
special properties.

This paper is organized as follows. In Section~\ref{prelims}, we 
recall necessary definitions and basic properties of Kohn 
Laplacians, the $\bar\partial$--Neumann problem, and Schr\"{o}dinger operators. Section~\ref{aux} 
contains some basic inequalities for $\complex^1$, including 
Lemma~\ref{magneticlemma}, which quantifies the key magnetic effect 
on which Theorem~\ref{main} is based. The 
example of Theorem~\ref{main} is constructed in 
Sections~\ref{setconstruction}  through \ref{muconstruction}. The 
verification that it possesses the desired properties is given in 
Section~\ref{endproof}.  
It will be apparent that Proposition~\ref{prop:smoothcase} follows from 
a simplification of the same construction, so we will not provide
a formal proof.
Theorem~\ref{thm:smoothcase} is proved in Section~\ref{section:smoothcase}.
Theorem~\ref{compactness} is then proved in 
Section~\ref{reduction} by reducing questions concerning property 
($P$) and compactness to semi-classical analysis of Schr\"{o}dinger 
operators.  
Finally, the reduction of Theorem~\ref{thm:dbarsmoothcase} to
Theorem~\ref{thm:smoothcase} is indicated in \S\ref{section:lastreduction}.


We are grateful to B.~Simon for clarifying for us the 
history of Kato's inequality, and of diamagnetic and paramagnetic inequalities.
The second author also thanks J.~J.~Kohn and  E.~Straube for encouragement and stimulating discussions.

\section{Preliminaries}\label{prelims}

\subsection{Kohn Laplacians and notions of compactness} Let $\Omega=\{(z, w):
\ \ |w|<e^{-\varphi(z)}, \ z\in\Omega_0\}$ be a complete Hartogs
domain in $\C^d\times C^1$ for some $d\ge 1$.
Assume that $b\Omega$ is both smooth and strictly
pseudoconvex in a neighborhood of $b\Omega\cap\{w=0\}$.  (This is
equivalent to the conditions that $\Omega_0$ has smooth boundary,
$\lim_{z\to b\Omega_0}\varphi(z)=\infty$, and there exists a
subdomain $\hat\Omega_0\subset\subset\Omega_0$ such that
$-e^{-2\varphi}$ has a smooth extension to a neighborhood
of $\overline{\Omega_0}$, whose complex Hessian is
strictly positive on $\overline{\Omega}_0\setminus\hat\Omega_0$.)
We will also assume that $\varphi$ is subharmonic on $\Omega_0$ 
(which is equivalent, under our other hypotheses, to pseudoconvexity
of $\Omega_0$) and $\partial\varphi/\partial z
\in L^2$ in any compact subset of $\Omega_0$,
in the sense of distributions.)
All this implies that there exists an open domain $\Omega'_0\Subset\Omega_0$
such that $\varphi$ is $C^\infty$ and strictly plurisubharmonic
in $\Omega_0\setminus\Omega'_0$ as well as in a neighborhood
of the boundary of $\Omega'_0$.

Let $M=b\Omega$. The portion of $M$ on which $z\in\Omega'_0$
will be parametrized by the projection $z\in\Omega'_0$
and thus identified with $\Omega'_0$.
For $1\le j\le d$
let $\bar L_j=\partial_{\bar z_j} -i\varphi_{\bar z_j}
\partial_\theta$
and let ${L_j}$ be its conjugate.  $L_j$ and
$\overline{L}_j$ may be considered as operators
defined in the sense of distributions
on $L^2(\Omega_0\times\T)$.
On $M$ one has formally the usual complex
of Cauchy-Riemann operators $\dbarb$,
mapping $(0,q)$ forms to $(0,q+1)$ forms.
We equip $M=b\Omega$ with a measure which has a nonvanishing
$C^\infty$ density with respect to the induced surface
measure wherever $z\notin\Omega'_0$, and which agrees
with Lebesgue measure in the coordinate $z\in \complex^d$
wherever $z$ lies in a neighborhood of the closure of $\Omega'_0$.
Note that surface measure, in contrast, carries a
factor related to $\nabla\varphi$ so that when $\nabla\varphi$
is merely square integrable, surface measure is not equivalent to
the measure which we have chosen.

Denote by $B^{0,q}$ the bundle of $(0,q)$ forms. Any section
can be expressed as $f = \sum_J f_J(z) \overline{dz}_J$,
and $\dbarb f = \sum_{j=1}^d \bar L_j f \overline{dz}_j\wedge
\overline{dz}_J$.
We choose a Hermitian metric for $B^{0,1}$ so that
$\{\overline{dz}_J\}$ form an orthonormal basis for $B^{0,q}$
at each point of $\Omega'_0$.

For $d\ge 2$ (that is, for domains in $\complex^3$), let
\[
Q_b(f, f)=\|\dbarb f\|_{L^2}^2 + \|\dbarb^* f\|_{L^2}^2
\]
for all sections $f$ of $B^{0,1}$ belonging to $C^1_0(\Omega'_0)$.
When $d\ge 2$, for smoothly bounded domains, this is
equivalent to the usual notion of compactness for $\square_b$,
the operator related to the closed sesquilinear form
$Q_b$ by $\langle \square_b f,\,f\rangle = Q_b(f,f)$;
this equivalence is a consequence of well-known estimates
since $M$ is strictly pseudoconvex where $z\notin\Omega'_0$.

For $d=1$ we define instead
\[
Q_b(f,f) = \|\bar L f\|_{L^2}^2 + \|L f\|_{L^2}^2,
\]
for scalar-valued $f\in C^1_0(\Omega'_0)$, 
where $\bar L = \bar L_1$.
\begin{definition} \label{compactnessdefn}
Let $\Omega\subset\complex^2$, that is, $d=1$.
We say that the Kohn Laplacian has compact resolvent in $L^2(b\Omega)$ if
the set of all $f\in C^1_0(\Omega'_0)$ for
which $Q_b(f,f)\le 1$ is precompact in $L^2(\Omega'_0)$.
\end{definition}

For $d=1$, that is for smooth domains in $\complex^2$,
an alternative notion of compactness is that
the set of all scalar-valued $f\in C^1$ which are orthogonal to the
$L^2$ nullspace of $\dbarb$ and satisfy $\|\dbarb f\|_{L^2(b\Omega)}
\le 1$ should be precompact in $L^2(b\Omega)$.
For smoothly bounded pseudoconvex domains in $\complex^2$, the range of $\dbarb$ in $L^2$
is known to be closed, and compactness in the sense of Definition~\ref{compactnessdefn}
would thus imply compactness in this alternative sense.

By formulating
compactness as in Definition~\ref{compactnessdefn},
we have avoided
discussing whether $\dbarb$ has closed range for
the class of nonsmooth Hartogs domains, Hilbert space, and Hermitian
structures investigated here; we have likewise sidestepped the question
of the relation between compactness for the boundary Kohn Laplacian,
and compactness for the $\bar\partial$-Neumann problem for the interior domain. 
In particular, the question of whether surface measure, or our alternative measure, 
is the relevant measure to place on the boundary has not been analyzed.
 
Matheos \cite{Matheos97} proved that
for arbitrary bounded pseudoconvex domains $\Omega\subset\complex^2$ with $C^\infty$ boundaries,
compactness holds in the $\dbar$-Neumann problem 
if and only if the following boundary estimate holds:  
For any $\varepsilon>0$, there exists $C_\varepsilon>0$ such that
\begin{equation}\label{2Dcompactnessestimate} 
\|u\|^2 \le 
\varepsilon(\|Lu\|^2+\|\overline{L}u\|^2)+C_\varepsilon
\|u\|^2_{-1} 
\end{equation} 
for all $u\in C^\infty(b\Omega)$.
This is equivalent to compactness in the sense of Definition~\ref{compactnessdefn}.



\subsection{The $\dbar$-Neumann Laplacian}  
Let $\Omega$ be a bounded domain in $\C^2$. Let $L^2_q(\Omega)$, $0\le q\le 2$,
be the space of $(0, q)$ forms with $L^2$-coefficients, equipped
with the standard Euclidean metric.

Let $\dbar_q\colon L^2_q(\Omega)\to L^2_{q+1}(\Omega)$ be defined
in the sense of distributions with $\dom(\dbar_q)=\{f\in L^2_q(\Omega):
\ \ \dbar_q f\in L^2_{q+1}(\Omega)\}$.  Let $\dbarstar_q$ be the adjoint of
$\dbar_{q-1}$.   Consider
\[
Q(u, v)=(\dbar_1 u, \dbar_1 v)+(\dbarstar_1 u, \dbarstar_1 v)
\]
with $\dom(Q)=\dom(\dbar_1)\cap\dom(\dbarstar_1)$.  It is easily to
see that $Q$ is a densely defined, non-negative, closed
sesquilinear form.  Therefore it uniquely defines a densely
defined, non-negative, self-adjoint operator $\square\colon
L^2_1(\Omega)\to L^2_1(\Omega)$ such that
$\dom(\square^{1/2})=\dom(Q)$ and $Q(u, v)=(\square u, v)$ for
$u\in \dom(\square)$ and $v\in \dom(Q)$. The operator $\square$ is
called the $\dbar$-Neumann Laplacian. It is said to have compact
resolvent if $(I+\square)^{-1}\colon L^2_1(\Omega)\to
L^2_1(\Omega)$ is compact.  This is in turn equivalent to the
following compactness estimate: For any $\eps>0$ there exists $C_\eps<\infty$ such that
\[
\|u\|^2\le \e Q(u, u)+ C_\e \|u\|^2_{-1}, \qquad \text{for all }u\in \dom(Q).
\]

It follows from the $L^2$-estimates of H\"{o}rmander  \cite{Hormander65}
for $\dbar$ that when $\Omega$ is pseudoconvex, $\square$ is 1-1 and onto, and
therefore has a bounded inverse $N$,  which is called the
$\dbar$-Neumann operator. In this case, $\square$ has compact
resolvent if and only if $N$ is compact.

\subsection{Schr\"odinger operators in $\complex^1$} 
 
Let $\psi$ be a subharmonic function defined in a bounded 
domain $\Omega_0\subset\complex^1$.  Assume that 
$\nabla\psi\in L^2(\Omega_0)$, in the sense of distributions. 
Let $D_\psi=(\p_x+i\psi_y, \p_y-i\psi_x)$.  Let 
\[ 
Q_\psi(u, v)=\langle D_\psi(u), D_\psi(v)\rangle 
+\langle u, \Delta\phi v\rangle 
\] 
be the closed, non-negative sesquilinear form on $L^2(\Omega_0)$ 
with core $C^\infty_0(\Omega_0)$.  This sesquilinear form uniquely 
defines a non-negative, self-adjoint, densely defined operator 
$S_\psi$ on $L^2(\Omega_0)$.  $S_\psi$ is the 
Schr\"{o}dinger operator with magnetic potential 
$A=(-\psi_y, \psi_x)=-\psi_y dx +\psi_x dy$, 
magnetic field $dA=\Delta\psi dx\wedge dy$, and electric 
potential $V=\Delta\psi$.  It is formally written as 
\begin{align*} 
S_\psi &=D_\psi^*\cdot D_\psi +\Delta \psi \\ &= 
-[(\p_x+i\psi_y)^2 + (\p_y -i\psi_x)^2] + \Delta\psi . 
\end{align*} 
Let $S^0_{\psi} = -\p_x^2-\p_y^2 + \Delta\psi$ be the Schr\"{o}dinger operator 
with the same electric potential but zero magnetic potential. 
The lowest eigenvalue of $S^0_\psi$ is given 
by 
\begin{definition}
\[ 
\lambda^{\rm e}_\psi=\inf\{ \|\nabla u\|^2+\|\sqrt{\Delta\psi} u\|^2; 
\ \ u\in C^\infty_0(\Omega_0), \|u\|=1\}. 
\] 
\end{definition}
\noindent
The lowest eigenvalue of $S_\psi$ is 
\begin{definition}
\begin{equation*}
\lambda^{\rm m}_\psi = \inf\{Q_\psi(u, u); 
\ \ u\in C^\infty_0 (\Omega_0), \|u\|=1\}. 
\end{equation*}
\end{definition}
\noindent
$\lambda^{\rm m}_\psi$ may alternatively be expressed as
\begin{align*} 
\lambda^{\rm m}_\psi &=
\inf\{4\|L_\psi (u)\|^2;
\ \ u\in C^\infty_0 (\Omega_0), \|u\|=1\} \\ 
                &=\inf\{4\|u_z e^\psi\|^2; 
\ \ u\in C^\infty_0 (\Omega_0), \|ue^\psi\|=1\} , 
\end{align*} 
where $L_\psi=-\p_z +\psi_z$. The last equality above follows 
from an easy substitution while the preceding
equality follows from the integration by parts formula: 
\begin{align*} 
\langle S_\psi(u), u\rangle =4\int_{\Omega_0} |L_\psi(u)|^2 
=\int_{\Omega_0} |D_\psi (u)|^2 + \int_{\Omega_0} \Delta\psi 
|u|^2. 
\end{align*} 
Another useful integration by parts formula is the following 
twistor formula. 
 
\begin{equation}\label{twistorformula} 
\int_{\Omega_0} a|L_\psi u|^2  =\int_{\Omega_0} 
\left((2a\psi_{z\bar z} -a_{z\bar z})|u|^2 + a|\bar L_\psi 
u|^2\right)  + 2\Re \int_{\Omega_0} u a_z \overline{L_\psi (u)} 
\end{equation} 
for any $a\in C^2(\overline{\Omega}_0)$.  Let $b\in 
C^2(\overline{\Omega}_0)$ and $b\le 0$.  Using the above formula 
with $a=1-e^b$ and applying the Schwarz inequality to the last 
term, we then obtain 
\begin{equation}\label{twistorinequality} 
\int_{\Omega_0} |L_\psi u|^2 \ge \int_{\Omega_0} b_{z\bar z} |u|^2 
e^b +\int_{\Omega_0} 2a\psi_{z\bar z} |u|^2  +\int_{\Omega_0} a 
|\bar L_\psi u|^2. 
\end{equation}

\section{Basic inequalities}  \label{aux} 
 
In this section, we collect several inequalities which will be 
used in the analysis. We start with the following well-known 
inequality of Kato (e.g.\  \cite{Kato72} and \cite{Simon79}), 
whose relevance to diamagnetism was observed by Simon \cite{Simon76}.
Integrals are taken with respect to Lebesgue measure 
on $\complex^1$, except where otherwise indicated. 
 
\begin{lemma} \label{Kato} 
Let $\psi$ be a real-valued function on a domain 
$\Omega_0\subset\complex^1$ such 
that $\nabla\psi\in L^2$, in the sense of distributions. 
Let $u\in C^1(\Omega_0)$. Then 
$\big|\nabla |u|(z)\big|\le |D_\psi u(z)|$ for {\it a.e.\ } 
$z\in \Omega_0$.  In particular, 
\begin{equation}  \label{gradientcontrol} 
\int_{\Omega_0} |D_\psi u|^2 \ge \int_{\Omega_0} \big|\nabla|u|\big|^2. 
\end{equation} 
\end{lemma} 
 
A short proof is provided for the reader's convenience. 
\begin{proof} 
$|u|$ is Lipschitz continuous, hence is differentiable almost 
everywhere. The $L^\infty$ function $\nabla|u|$ thus defined 
equals the gradient in the distribution sense, and $\big|\nabla|u| \big| 
\le |\nabla u|$ a.e. 
 
At points where $u$ vanishes, the magnetic gradient  
equals the ordinary gradient, so the conclusion holds. 
In the open set where $u\ne 0$, one can locally write 
$u(z) = r(z)e^{i\theta(z)}$ with $r,\theta\in C^1$. 
Then $\nabla|u|=\nabla r$, while 
$(\partial_x+i\psi_y)re^{i\theta} 
= (r_x+i\psi_y r)e^{i\theta}$ has magnitude 
$(|r_x|^2 + |\psi_y r|^2)^{1/2} 
\ge |r_x|$. Bounding $(\partial_y -i\psi_x)re^{i\theta}$ 
in the same way leads to the desired inequality. 
\end{proof} 
 
 
For any $x\in\reals$ define 
\begin{equation} 
\|x\|_* = \distance(x,\integers). 
\end{equation} 
 
The next lemma indicates one situation in which 
the magnetic gradient is relatively powerful; 
in fact it will be the key ingredient in our proof that 
$\lambda^{\rm m}_{n\varphi}\to\infty$. 
The result is also not original; for much more general
results of the same type see \cite{laptevweidl} and \cite{balinski}.
 
\begin{lemma}  \label{magneticlemma} 
If $\Delta\psi\equiv 0$ in 
$\scripta=\{z: r<|z|<R\}$, then for any $u\in C^1$, 
\begin{equation}  \label{annularlowerbound} 
\int_\scripta |D_\psi u|^2 
\ge \|w(\psi)\|_*^2  \int_\scripta |z|^{-2} |u(z)|^2, 
\end{equation} 
where the winding number $w(\psi)$ is given by 
\begin{equation}  \label{windingnumberdefn} 
w(\psi) = \frac1{2\pi}\int_{|z|=\rho} -\psi_y dx +\psi_x dy 
\end{equation} 
for any $\rho\in (r, R)$. 
More precisely, 
\begin{equation}  \label{circlelowerbound} 
\int_0^{2\pi} |D_\psi u(\rho e^{i\theta})|^2\,d\theta 
\ge \rho^{-2}\|w(\psi)\|_*^2  \int_0^{2\pi} 
|u(\rho e^{i\theta})|^2\,d\theta 
\end{equation} 
for any $\rho\in(r,R)$. 
\end{lemma} 

This expresses one instance of the Aharonov-Bohm phenomenon. The magnetic
field $\Delta\psi$ vanishes identically in $\scripta$, yet if $w(\psi)\ne 0$
then (roughly speaking) a quantum particle confined to $\scripta$  
and governed by the Hamiltonian $D_\psi^* D_\psi$ experiences a measurable
effect from the magnetic potential. For a semiclassical analysis
of this effect in certain cases see \cite{helfferaharonovbohm}.
 
By Stokes' theorem together with the assumption $\Delta\psi=0$, 
the integral \eqref{windingnumberdefn} defining the winding number 
is independent of $\rho\in(r,R)$. 
If $\psi$ extends to a $C^2$ function in the disk $|z|<R$, harmonic 
where $|z|>r$, then there is the alternative expression 
\begin{equation}  \label{altwindingnumberexpression} 
w(\psi) = (2\pi)^{-1} \int\Delta\psi\,dx\,dy. 
\end{equation} 
 
\begin{proof} 
Using the polar coordinates $z=re^{i\theta}$, a straightforward 
calculation gives 
\begin{equation}\label{polarcoordinates} 
|D_\psi u|^2=|u_r +ir\rp  \psi_\theta 
u|^2 +|r\rp u_\theta -i\psi_r u|^2. 
\end{equation} 
It suffices to prove \eqref{circlelowerbound}, which directly 
implies \eqref{annularlowerbound}. 
Let $\tilde\psi$ be the harmonic conjugate of 
$\psi -w(\psi)\log |z|$ on $\scripta$.  Let $v=ue^{-i\tilde\psi}$.  Then 
$|D_\psi u|=|D_{w(\psi)\log |z|} v|$, and $|v|\equiv |u|$. 
Write 
\[ 
v(\rho e^{i\theta})=\sum_{k=-\infty}^{\infty} \hat v(k,\rho) e^{ik\theta} 
\] 
where this expression defines the Fourier coefficients $\hat v$. 
It follows from~\eqref{polarcoordinates} that 
\begin{align*} 
\int_0^{2\pi} |D_{w(\psi)\log |z|} v(\rho e^{i\theta})|^2 d\theta 
&\ge \rho^{-2}\int_0^{2\pi} |v_\theta -iw(\psi) v|^2  d\theta\\ 
&= \rho^{-2} \int_0^{2\pi} |\sum_{k=-\infty}^{\infty} \hat 
v(k,\rho) 
(k-w(\psi))e^{ik\theta}|^2 d\theta \\ 
&\ge \rho^{-2}\|w(\psi)\|^2_{*} \int_0^{2\pi} |v(\rho 
e^{i\theta})|^2 d\theta. 
\end{align*} 
The lemma then follows. 
\end{proof} 
 
A more general result holds, although only the special 
case formulated in Lemma~\ref{magneticlemma} will be needed in our analysis. 
 
\begin{lemma}  \label{new3.2} 
Let $\Gamma$ be a rectifiable Jordan curve of length $\rho$, 
parametrized by arclength $s\in [0,\rho]$. 
Let $h$ be a real-valued function on $\Gamma$, regarded 
as a function of $s$. Let $L$ be the first-order differential 
operator $\frac{d}{ds}+ih$, acting on the space 
$L^2$ of {\em periodic} functions on $[0,\rho]$. 
Define the winding number $w=(2\pi)\rp \int_0^\rho h(s)\,ds$. 
Then for any periodic test function $u\in C^1([0,\rho])$, 
\begin{equation} 
\|Lu\|_{L^2} \ge 4\|w\|_*\rho^{-1}\|u\|_{L^2}. 
\end{equation} 
\end{lemma} 
 
\begin{proof} 
Writing $H(s) = \int_0^s h$ and $Lu=f$, we have $L = e^{-iH}\frac{d}{ds}e^{iH}$ 
so $\frac{d}{ds}(e^{iH}u) = e^{iH}f$, 
whence  
\begin{equation} 
e^{iH(s)}u(s)= e^{iH(0)}u(0) + \int_0^s e^{iH}f. 
\end{equation} 
Therefore 
\begin{equation} 
|e^{i(H(s)-H(0))} u(s)-u(0)| 
\le \int_0^s|f|. 
\end{equation} 
Applying the same reasoning to the interval $[s,\rho]$ gives 
\begin{equation} 
|u(\rho)-e^{i(H(s)-H(\rho))} u(s)| 
\le \int_s^\rho|f|. 
\end{equation} 
By the triangle inequality and the periodicity assumption $u(\rho)=u(0)$, 
this implies 
\begin{equation} 
|u(s)|\cdot |e^{i(H(s)-H(0))}- e^{i(H(s)-H(\rho))}|\le \int_0^\rho|f|, 
\end{equation} 
which is equivalent to 
\begin{equation} 
|u|_{L^\infty}\cdot |e^{iH(\rho)}-e^{iH(0)}|\le\|f\|_{L^1}. 
\end{equation} 
Now $|e^{iH(\rho)}-e^{iH(0)}|\ge 4\|w\|_*$, 
and applying Cauchy-Schwarz twice gives two factors of $\rho^{1/2}$. 
\end{proof} 
 
This implies \eqref{circlelowerbound},  
except for a constant factor in the inequality, 
by taking $L$ to be the component of the magnetic gradient 
tangent to $\Gamma$.  
For a general Jordan curve $\Gamma$, the winding number $w$ 
which appears in Lemma~\ref{new3.2} equals $\pi\rp\int_{\scriptr} 
\Delta\psi$, where $\scriptr$ is the region enclosed by $\Gamma$.

We will also need the following Poincar\'{e}-type inequalities. 
Denote by $B(z, r)$ the disk centered at $z$ with radius $r$. 
 
\begin{lemma}\label{poincare}   If $u\in C^0(\ov{B(0, R)}) 
\cap W^1(B(0, R))$, then 
\begin{equation}\label{poincare1} 
\int_{B(0, R)} |u|^2 \le R^2\left(2\int_{B(0, R)} \big|\nabla 
|u|\big|^2 +\int_0^{2\pi} |u(Re^{i\theta})|^2 d\theta\right). 
\end{equation} 
Let $\scripta =\{r<|z|<R\}$. 
If $u\in C^0(\ov{\scripta})\cap W^1(\scripta)$ then 
\begin{equation} \label{poincare2} 
\int_\scripta |u|^2 \le (R^2-r^2)\log(R/r) \int_\scripta 
\big|\nabla |u|\big|^2 + (R^2-r^2)\int_0^{2\pi} 
|u(re^{i\theta})|^2\,d\theta. 
\end{equation} 
\end{lemma} 
 
\begin{proof}  $|u|$ likewise belongs to $C^0\cap W^1$. 
Using Friederichs mollifiers permits us to assume that $|u|\in C^1$. 
Therefore, without loss of generality, we assume that $u\ge 0$ 
and $u\in C^1$. 
 
To prove \eqref{poincare2} 
we work in polar coordinates, and exploit only the radial component of the gradient. 
It thus suffices to show that 
\begin{equation*} 
\int_r^R |f(\rho)|^2\rho\,d\rho \le \log(R/r)(R^2-r^2)\int_r^R 
|f'(\rho)|^2\,\rho\,d\rho + (R^2-r^2)|f(r)|^2 
\end{equation*} 
for any $f\in C^1(\reals)$. 
Express $f(\rho) = f(r)+\int_r^\rho f'(t)\,dt$ 
and note that $\int_r^R|f(r)|^2\,\rho\,d\rho= \tfrac12 (R^2-r^2)|f(r)|^2$. 
The other term is 
\begin{multline*} 
\le 2\int_r^R \big(\int_r^\rho |f'(t)|\,dt \big)^2 \rho\,d\rho 
\\ 
\le 2\int_r^R \big(\int_r^\rho |f'(t)|^2t\,dt \big) \big(\int_r^\rho t\rp\,dt \big) 
\rho\,d\rho 
\le \log(R/r) (R^2-r^2)\int_r^R |f'|^2\,t\,dt, 
\end{multline*} 
as claimed. 
The proof of \eqref{poincare1} is similar and is left to the reader. 
\end{proof} 
 
\section{Construction of thick sets in $\complex^1$} 
\label{setconstruction} 
 
In this section we explicitly construct sets in $\C^1$ that have 
empty Euclidean interior and non-empty fine interiors.  This 
construction will later be used in the proof of 
Theorem~\ref{main}. 
 
Let $B$ be an integer greater than 2.  We always assume that $B$
is chosen to be sufficiently large for various inequalities encountered below
to be valid. For any positive integer 
$k$, let $\e_k=B^{-k}$ and let $\Lambda_k=B^{-k}\cdot(\ZZ +i\ZZ)$ 
be the set of lattice points.  Let $\rho_k$ be a positive number 
of the form
\begin{equation}  \label{rhoconstraint}
\rho_k = e^{-\sigma_k k B^{2k}}
\ \text{ where } \sum_k \frac1{k\sigma_k}<\infty  \text{ and } \sigma_k\ge 1.
\end{equation}
In particular, $\rho_k$ is much smaller than any
power of $\e_k$ for large $k$.

Let $\Omega_0$ be the 
unit disk.  We choose $\{z^k_j\}_{j=1}^{m_k}\in 
\Lambda_k\cap\Omega_0$ by induction on $k$ as follows.  For $k=1$, 
$\{z^1_j\}_{j=1}^{m_1}$ are chosen to be all points in 
$\Lambda_1\cap\Omega_0$ such that $\ov{B(z^1_j, 
\e_1)}\subset\Omega_0$.  Suppose the points $z^l_j$ have been 
chosen for $l\le k-1$. Then $\{z^k_j\}_{j=1}^{m_k}$ are chosen to 
be all those points in $\Lambda_k\cap\Omega_0$ such that 
$\ov{B(z^k_j, \e_k)} \subset\Omega_0$ and 
\[ 
\distance(z^k_j, \cup_{l=1}^{k-1}\cup_{j=1}^{m_l} D^l_j) > \e_k 
\] 
where $D^l_j=B(z^l_j, \rho_l)$. Note that $m_k\le 4\e^{-2}_k$.  Let 
\[ 
\Omega_k = \Omega_0\backslash\cup_{l=1}^k\cup_{j=1}^{m_l} 
\ov{D^l_j} \quad \text{and}\quad W = 
\Omega_0\backslash\cup_{l=1}^\infty \cup_{j=1}^{m_l} 
\ov{D^l_j}=\cap_{k=1}^\infty \Omega_k. 
\] 
It is evident that the set $\Omega$ thus constructed has empty (Euclidean) 
interior. Moreover, $\Omega$ has positive Lebesgue measure because 
\[ 
\sum_{k=1}^\infty\sum_{j=1}^{m_k} |D^k_j|\le4\pi\sum_{k=1}^\infty 
\rho_k^2\e_k^{-2}\le 4\pi\sum_{k=1}^\infty \e_k^2 \ll |\Omega_0|
\] 
provided that $B$ is chosen sufficiently large.
\begin{figure}[htbp] \centering
 \ \psfig{figure=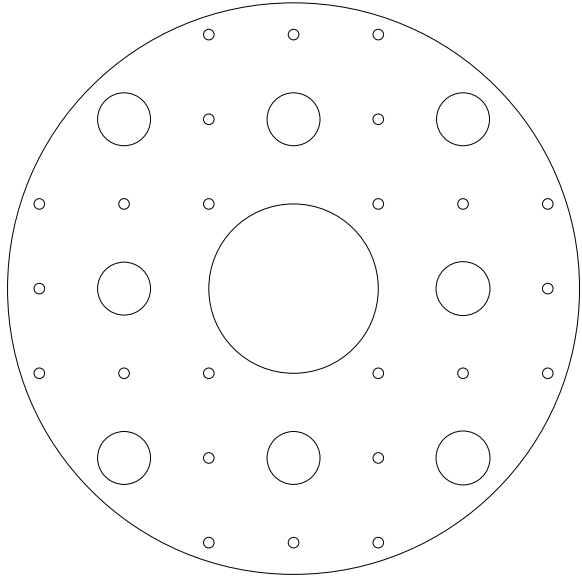,width=3in,height=3in}
 \caption{
The domain $\Omega_0$ and disks $D^k_j$ of generations $k=1,2,3$.
Their centers $z^k_j$ lie on a lattice of scale $B^{-k}$,
whereas their radii $\rho_k$ approach $0$ at a doubly exponential rate.}
 \end{figure}
 
\begin{lemma}\label{counting} 
When $B$ is sufficiently large, $m_k\ge 
\e_k^{-2}/4$. Furthermore, for any $z\in\Omega_k$ satisfying 
$\distance(z,\partial\Omega_0)\ge\e_k^{1/2}$, there exist $\ge 
\e_k^{-1}$ indices $j$ such that $|z^k_j -z|< 4\e_k^{1/2}$. 
\end{lemma} 
 
\begin{proof} We use the simple fact that for any $r\ge 1$, the number of integer 
lattice points in a (closed or open) disk of radius $r$ is bounded 
between $r^2$ and $5r^2$.  It follows that the number of points of 
$\Lambda_k\cap\Omega_0$ that are not elements of $\{z^k_j\}$ is no 
more than 
\[ 
\sum_{1\le l<k} 4\e_l^{-2}\cdot 
5\big(\frac{\e_k+\rho_l}{\e_k}\big)^2 \le 40\sum_{1\le l<k} B^{2l} 
\big( 1 + B^{2k-4l} \big) \le C \e_k^{-2} B^{-2}. 
\] 
Since the cardinality of $\Lambda_k\cap \{|z|<1-\e_k\}$ is $\ge 
\e_k^{-2}/2$, we have $m_k\ge (1/2 -C B^{-2}) 
\e_k^{-2}>\e_k^{-2}/4$. 
 
We now prove the second statement. Indeed, the above reasoning 
still applies, unless $B(z,4\e_k^{1/2})$ meets some $D^l_i$ with 
$\e_l > 16\e_k^{1/2}$. In that case there can be only one such $l$ 
and only one such $i$.  Otherwise, consider the largest such $l$; 
by construction, the distance from $D^l_i$ to any $D^m_j$ with 
$m\le l$ is $\ge \e_l-2\rho_l\ge \e_l-2\e^2_l$. Thus $\e_l-2\e_l^2 
\le 8 \e_k^{1/2}<\e_l/2$, which is impossible.  Since $B(z, 
4\e_k^{1/2})$ meets only one $D^l_i$ with $\e_l>16 \e_k^{1/2}$, it 
must contain a disk $B(z', 2\e_k^{1/2})$ which does not intersect 
any $D^l_i$ with $\e_l>16 \e_k^{1/2}$.  As in the preceding 
paragraph, the number of points in $\Lambda_k\cap B(z', 
\frac3{2}\e_k^{1/2})$ that are not elements of $z^k_j$ is no more 
than 
\[ 
\sum_{\begin{subarray}{l} 
1\le l<k \\ 
\e_l\le 16\e_k^{1/2} \end{subarray}} 
C\big(\frac{2\e_k^{1/2}}{\e_l}\big)^2\cdot 
5\big(\frac{\e_k+\rho_l}{\e_k}\big)^2 \le C\e_k \sum_{1\le 
l<k}\e_l^{-2}\cdot \big(\frac{\e_k+\rho_l}{\e_k}\big)^2\le C
\e_k^{-1} B^{-2}. 
\] 
Since the cardinality of $\Lambda_k\cap B(z', \tfrac32\e_k^{1/2})$ 
is $\ge \frac9{4}\e_k^{-1}$, the number of points $z^k_j$ in $B(z, 
4\e_k^{1/2})$ is $\ge (\tfrac94 -C{B^{-2}})\e_k^{-1}\ge \e_k^{-1}$. 
\end{proof} 
 
\begin{lemma}\label{fine} 
If $\sum_k (k\sigma_k)^{-1}<\infty$
then there exist constants $c,C\in(0,\infty)$ and functions  
$F_k\in W^1_0(\Omega_k)$  
such that for all $k$, $\|F_k\|_{L^2(\Omega_0)}\ge c$, 
$\|\nabla F_k\|_{L^2(\Omega_0)}\le C$, and $\|F_k\|_{L^\infty}\le C$.   
\end{lemma} 
\begin{proof} 
Define 
\[ 
f^l_j(z) = \frac{\log(|z-z^l_j|/\rho_l)} 
{\log({\e_l^2}/4\rho_l)} 
\] 
if $\rho_l\le|z-z^l_j|\le {\e_l^2}/4$, and $f^l_j(z)\equiv 1$ 
for $|z-z^l_j|\ge {\e_l^2}/4$.  Note that $0\le f^l_j(z)\le 1$ for all $z$. Let 
\[ 
F_k(z)=(1-|z|^2)\min_{1\le l\le k}\min_{1\le j\le m_l} f^l_j(z) . 
\] 
If $z\in B(0, 1/2)$ and $|z-z^l_j|\ge {\e_l^2}/4$ for 
all $1\le l\le k$ and all $1\le j\le m_l$, then $F_k(z)=1-|z|^2\ge 3/4$. 
Moreover $0\le F_k(z)\le 1$ for all $z\in\Omega_0$.
Since 
\[ 
\sum_{k=1}^\infty\sum_{j=1}^{m_k} |B(z^k_j, {\e_k^2}/4)| 
\le C\pi\sum_k \e_k^{-2} \e_k^4, 
\] 
we have $\|F_k\|^2 \ge \pi/2$ provided that $B$ is chosen to
be sufficiently large.
 
In estimating $\|\nabla F_k\|$ 
from above, we may disregard the harmless factor of $1-|z|^2$. 
$\Omega_k$ can be partitioned into finitely many pairwise disjoint 
subregions such that $F_k$ is identically equal to some 
$f_j^l$ on each subregion, and such that the supports of $\nabla f^l_j$ 
are mutually disjoint for a fixed $l$.  Moreover, since $F_k$ is continuous, 
$\|\nabla (1-|z|^2)^{-1} F_k\|^2$ equals the sum of the squares of the $L^2$ norms of 
its gradients over all these subregions. Therefore since 
$F_k$ is bounded above in the supremum norm uniformly in $k$, 
\[ 
\|\nabla F_k\|^2\le C + \sum_{l=1}^k \sum_{j=1}^{m_l}\int_{\Omega_l} 
|\nabla f^l_j|^2. 
\] 
Since 
\[ 
\int_{\Omega_l} |\nabla f^l_j|^2= 2\pi 
\int_{\rho_l}^{{\e_l^2}/4} 
r^{-2}[\log({\e_l^2}/4\rho_l)]^{-2} \,r\,dr 
=2\pi[\log({\e_l^2}/4\rho_l)]^{-1}\le 
4\pi(\log(1/\rho_l))^{-1}, 
\] 
we have 
\[ 
\|\nabla F_k\|^2\le C\sum_{l=1}^k 
\e_l^{-2}\cdot  (\log(1/\rho_l))^{-1}
= C\sum_{l=1}^k B^{2l}\cdot B^{-2l}(l\sigma_l)^{-1}
= C\sum_{l=1}^k(l\sigma_l)^{-1}.
\]
\end{proof} 

 
\begin{remark}  The fine topology is the smallest 
topology on $\C$ with respect to which all subharmonic functions are 
continuous. We refer the reader to \cite{Helms69}, Chapter 10, for 
an elementary treatise on the fine topology.  It follows from 
\cite{Fuglede99} that the existence of functions $F_k$ satisfying
the conclusions of Lemma~\ref{fine} is equivalent to $W$ having nonempty fine 
interior. 
\end{remark} 
 
\section{Construction of the subharmonic function $\varphi$} 
\label{varphiconstruction} 
 
We follow the construction in the preceding section, taking 
$\rho_k =\exp(-k\sigma_k B^{2k})$.
Let $h\in C^\infty_0(\C)$ be the radially 
symmetric function defined by $h(t)=c_0 e^{-1/(1-t)}$ for $0\le 
t<1$ and $h(t)=0$ for $t\ge 1$, where the constant $c_0$ is chosen 
so that $\int_{\C} h=1$. 
For $k\in\naturals$ and $1\le j\le m_k$, define
\[ 
\psi^k_j = 2\pi \mu^k_j\rho_k^{-2} h((z-z^k_j)/\rho_k), 
\] 
where $\mu^k_j > 0$ 
satisfies 
\begin{equation} \label{muconstraint}
\mu_j^k \le\nu_k\e_k^2
\end{equation}
and the factors $\nu_k$ are chosen to satisfy
\begin{equation} \label{nusigmaconstraint}
\sum_k (1+\sigma_k)k\nu_k<\infty.
\end{equation}
Define also
\[ 
\varphi^k_j = \frac{1}{2\pi}\log|z|*\psi^k_j. 
\] 
Thus $\Delta\varphi^k_j=\psi^k_j$. 
Set 
\[ 
\varphi = \sum_k\sum_j\varphi^k_j; 
\] 
the next lemma guarantees convergence of this sum. 
 
\begin{lemma} \label{lemma:continuityofphi} 
If the sequences $\sigma_k,\nu_k$ satisfy \eqref{nusigmaconstraint},
then $\varphi$ is a subharmonic function on 
$\complex$, $\Delta\varphi$ is lower semicontinuous,
$\varphi\in C^0$, and $\nabla\varphi\in L^2(\Omega_0)$. 
\end{lemma} 
 
\begin{proof} Let $\psi=\sum_k\sum_j \psi_j^k$.  Then 
\[ 
\|\psi\|_{L^1}=\sum_k\sum_j \|\psi_j^k\|_{L^1} 
=\sum_k\sum_j \mu_j^k \lesssim \sum_k B^{2k}\nu_k\e_k^2
=\sum_k \nu_k<\infty, 
\] 
since $\sum_k k\nu_k$ converges by \eqref{nusigmaconstraint}.

Write 
\[ 
\varphi^k_j(z)=\mu^k_j \int_{|t|<1} 
\log |\rho_k t +z^k_j -z| h(t) . 
\] 
If $|z-z^k_j|>2\rho_k$, then 
\[ 
|\varphi^k_j (z)| \le 
C\mu^k_j \log(1/\rho_k). 
\] 
If $|z-z^k_j|\le 2\rho_k$, then likewise 
\begin{equation*} 
|\varphi^k_j (z)| 
\le 
C\mu^k_j 
\big( \log(1/\rho_k) 
+\int_{|t|<1} |\log(t+(z^k_j-z)/\rho_k)|\big) 
\le C\mu^k_j \log(1/\rho_k). 
\end{equation*} 
Moreover, whenever  
$z\notin B(z^k_j, \tfrac14\e_k)$, there is an improved
bound $|\varphi^k_j(z)|\lesssim \mu^k_j\log(1/\e_k)\lesssim \nu_k\e_k^2 k$. 
For any $z,k$ there exists at most one index $j$ for which
$|z-z^k_j|\le \e_k/4$.
It follows that for any $z\in\complex$, uniformly for any $k$,
\begin{equation*} 
\sum_j |\varphi^k_j (z) |
\lesssim \nu_k\e_k^2 \sigma_k kB^{2k}
+ \sum_{j: |z-z^k_j|\ge \e_k/4} \nu_k\e_k^2 k
\lesssim \nu_k\e_k^2 \sigma_k kB^{2k} 
+  B^{2k} \nu_k\e_k^2 k,
\end{equation*}
and consequently
\begin{equation*} 
\sum_k\sum_j |\varphi^k_j (z) |
\lesssim \sum_k (1+\sigma_k)k\nu_k<\infty
\end{equation*} 
since $\e_k = B^{-k}$.
Therefore the series defining $\varphi$ is uniformly convergent to a continuous function. 
Since each $\varphi^k_j$ is subharmonic, so is $\varphi$. 
 
We now estimate the $L^2$ norm of $\varphi$.
Write 
\[ 
\frac{\p\varphi^k_j}{\p z} =\tfrac12 \mu^k_j\int_{|t|<1} 
\frac{h(t)}{z-z^k_j-\rho_k t} . 
\] 
If $|z-z^k_j|>2\rho_k$, then 
\[ 
\left|\frac{\p\varphi^k_j}{\p z}\right| \lesssim 
\mu^k_j
|z-z^k_j|^{-1}. 
\] 
If $|z-z^k_j|\le 2\rho_k$, then 
\begin{equation*} 
\left|\frac{\p\varphi^k_j}{\p z}\right| \le 
\tfrac12 \mu^k_j
\int_{|t|<1} \frac{h(t)}{\rho_k 
|t-(z-z^k_j)/\rho_k|} 
\lesssim
\mu^k_j\rho_k^\rp
\int_{|t|<3}\frac{1}{|t|} \lesssim \mu^k_j \rho_k^\rp . 
\end{equation*} 
Therefore 
\begin{equation*} 
\|\nabla\varphi^k_j\|_{L^2(B(0, 1))} 
\le C\mu^k_j +C\mu^k_j \big(\int_{\rho_k}^{1} r^{-2}\,rdr \big)^{1/2}
\le C\mu^k_j(\log(1/\rho_k))^{1/2}.
\end{equation*} 
As in the estimation of the supremum norm of $\varphi$, there is a stronger
inequality
\[
\|\nabla\varphi^k_j\|_{L^2\{|z-z^k_j|\ge \e_k/2\}}
\lesssim \mu^k_j (\log(1/\e_k))^{1/2}\lesssim \mu^k_j k^{1/2}.
\]

Since there are at most $CB^{2k}$ indices $j$ for each $k$, 
\begin{multline*} 
\|\nabla\varphi\|_{L^2}
\lesssim 
\sum_k \Big(B^{2k} (\nu_k\e_k^2)^2 \log(1/\rho_k) \Big)^{1/2}
+ \sum_k\sum_j \nu_k\e_k^2 k^{1/2}
\\
\lesssim\sum_k
\big( B^{2k} \nu_k^2\e_k^4 \sigma_k kB^{2k} \big)^{1/2}
+ \sum_k B^{2k} \nu_k\e_k^2 k^{1/2}
\lesssim \sum_k \nu_k k^{1/2}\sigma_k^{1/2} 
+ \sum_k \nu_k k^{1/2}.
\end{multline*}
The hypothesis \eqref{nusigmaconstraint} guarantees convergence of these sums.
\end{proof}

\begin{remark} 
In order for $\Delta\varphi$ to be $C^\infty$, or even H\"older continuous,
it is necessary that the far more restrictive condition $\mu^k_j = O(\rho_k^\alpha)$ 
for some $\alpha>2$ be satisfied. 
\end{remark} 
 
\section{The coefficients $\mu^k_j$} 
\label{muconstruction} 
 
Let $B$ be any fixed positive integer, sufficiently large so that the hypotheses
of Lemmas \ref{counting} and \ref{fine} are satisfied,
and recall that $\e_k=B^{-k}$ and $\rho_k=\exp(-k\sigma_k B^{2k})$. 
The basic strategy in the 
proof of Theorem~\ref{main} is to combine \eqref{circlelowerbound} 
with \eqref{gradientcontrol}, using the former to gain a strong 
bound over many circles, and the gradient estimate from 
\eqref{gradientcontrol} to then gain control over the remainder of 
$\Omega_0$. For each large $n$, we want to find lots of (disjoint) 
circles $\Gamma$, for which \eqref{circlelowerbound}, applied to 
$\psi=n\varphi$, gives a strong lower bound 
on $\int_\Gamma |u|^2$. (These circles can have different centers.) 
Then we use \eqref{gradientcontrol} on the complement of the union 
of all these circles, with 
\eqref{gradientcontrol} giving us good control on the boundary of 
the complement. 
 
The factor of $\rho^{-2}=|z|^{-2}$ in \eqref{circlelowerbound} 
is important, since it tends to make $\int|u|^2$ much smaller 
than $\int |D_\psi u|^2$, 
provided the circle has small radius. It can also be 
used to gain satisfactory control of the boundary terms in 
\eqref{poincare1} and \eqref{poincare2}, if for instance the annulus $r/2 \le |z|\le r$ 
is one on which \eqref{annularlowerbound} gives a good bound 
on $u$. On the other hand, 
we lose something in applying \eqref{poincare2} to annuli 
for which $\log(R/r)$ is too large, relative to $R^2$. 
Thus for each large $n$, we want to have a large number of 
such good circles, and we want them to be fairly densely distributed 
in the sense that for each $z$ and $n$ there is such a circle within 
distance $b_k$ of $z$, where $k=k(n)$ and $b_k\to 0$ 
at some rate to be specified. 
 
We use the following setup. There will be a sequence of positive 
integers $N_k$ converging rapidly to $+\infty$. 
To each $n\in[N_k,N_{k+1})$ 
we will associate a family $\scriptf(n)$ of disks $D^k_j$ with the 
following properties: 
\newline 
(1) For any $z\in\Omega_0$ there exists $D^k_j\in\scriptf(n)$ 
satisfying $\distance(z,D^k_j)\le C\e_k^{1/2}$. 
\newline 
(2) For any $D^k_j\in\scriptf(n)$, $\|n\mu^k_j\|_*\ge \tfrac14$, 
where $\|x\|_* =\distance(x,\integers)$. 
A given disk $D^k_j$ is permitted to belong 
to $\scriptf(n)$ for many different values of $n$. 
 
Define 
\begin{equation} \label{Nsubkdefn} 
N_k= 2^{B^{k-1}} 
\end{equation} 
For any sufficiently large $k$, we construct $\{\scriptf(n): 
N_k\le n<N_{k+1}\}$, and $\{\mu^k_j \}$ as follows.  We first 
cover $\Omega_{k-1}$ by $\sim \e_k^{-1}$ disks centered in 
$\Omega_k$ with radii $8\e_k^{1/2}$. We may arrange these covering 
disks so that each $D^k_j$ belongs to at least one covering disk 
and the shrinking by half of each covering disk is disjoint from 
the other covering disks.  By Lemma~\ref{counting}, we can then 
partition the disks $D^k_j$ into $\sim \e_k^{-1}$ subfamilies, 
each of cardinality $\ge \e_k^{-1}$, so that for each subfamily, 
all of its member disks are contained in a common covering disk. 
 
The number $\mu^k_j$ are chosen as follows, to ensure the existence of 
many disks with favorable winding numbers $n\mu_j^k$ for each integer 
$n\in[N_k,N_{k+1})$. 
Consider first $n=N_k$. Choose one disk from each subfamily, 
and let $\scriptf(n)$ be the set of all disks thus chosen. 
For each disk $D^k_j\in\scriptf(n)$, 
define the weight $\mu^k_j$ by 
\begin{equation} 
n\mu^k_j=\tfrac14. 
\end{equation} 
Then 
$\|n\mu^k_j\|_*=\tfrac14$; moreover, 
$\|m\mu^k_j\|_* \ge 1/4$ for all $m\in [N_k, 2N_k)$. 
For each $m\in[N_k,2N_k)$, set $\scriptf(m)=\scriptf(N_k)$. 
 
Next consider $n = 2N_k$, and repeat the procedure: let 
$\scriptf(n)$ be a collection consisting of one disk $D^k_j$ from 
each subfamily, not previously chosen. Define $\mu^k_j$ by 
$n\mu^k_j=\tfrac14$. Then $\|n\mu^k_j\|_*\ge\tfrac14$ for all 
$n\in[2N_k,4N_k)$. For $m\in [2N_k,4N_k)$, let 
$\scriptf(m)=\scriptf(2N_k)$. The next iteration begins with 
$n=4N_k$, and so on. Repeat the procedure until every integer 
$n\in[N_k,N_{k+1})$ has been considered, and $\scriptf(n)$ and 
associated coefficients $\mu^k_j$ defined. There are sufficiently 
many disks $D^k_j$ to allow this because $N_k\cdot 
2^{\e_k^{-1}}>N_{k+1}$. Any disks $D^k_j$ not in 
$\cup_{n=N_k}^{N_{k+1}}\scriptf(n)$ play a lesser role in the 
analysis; we set $\mu^k_j = N_{k+1}^{-1}$ for those although any 
sufficiently small strictly positive quantity would suffice. 
 
Thus
\[ 
\mu^k_j\le N_k^{-1} = 2^{-B^{k-1}}
\] 
for all sufficiently large $k$. 
On the other hand, we have already imposed the constraint \eqref{muconstraint}
$\mu^k_j\le\nu_k\e_k^2 = \nu_k B^{-2k}$,
with which the above construction is consistent if 
\[
\nu_k\ge B^{2k} 2^{-B^{k-1}}.
\]
In order to apply Lemma~\ref{lemma:continuityofphi} to conclude that
$\varphi \in C^0$ and $\nabla\varphi \in L^2$, 
we also need the constraints \eqref{rhoconstraint} $\sum k\rp \sigma_k^{-1}<\infty$
and \eqref{nusigmaconstraint} $\sum_k (1+\sigma_k)k\nu_k<\infty$.
All these are mutually compatible. Indeed if we fix any $\e>0$ and set
$\sigma_k = k^\eps$ and $\nu_k = k^{-2-2\eps}$ for large $k$, 
then $\nu_k\ge B^{2k} 2^{-B^{k-1}}$ with some room to spare.

The conclusions of this section are summarized in the following lemma.
\begin{lemma}
Suppose that $\nu_m\ge B^{2m} 2^{-B^{m-1}}$ for all sufficiently large $m$.
Then there exist coefficients $0\le\mu^m_j$ satisfying \eqref{muconstraint}
such that for each sufficiently large positive integer $n$
there exist an index $k=k_n$ and a collection $J_n$ of indices $j$ such that
$k_n\to\infty$ as $n\to\infty$,
and such that for each point $z\in \Omega_{k-1}$ there exists at least one $j\in J_n$
such that $|z-z^{k}_j|\le C\e_k^{1/2}$ and $\distance(\pi\rp n\int_{D^k_j}\Delta\varphi,
\integers)\ge\tfrac14$.

Moreover it is possible to choose a sequence $(\nu_m)$ and an associated sequences $(\sigma_m)$
such that
\eqref{rhoconstraint} and \eqref{nusigmaconstraint} are also satisfied.
Therefore $\varphi$ is subharmonic, $\varphi \in C^0$ and $\nabla\varphi \in L^2$,
and there exist functions $F_k$ satisfying the conclusions of Lemma~\ref{fine}.
\end{lemma}
 
\section{Proof of Theorem~\ref{main}} 
\label{endproof} 
 
We now proceed to prove that $\lambda^{\rm m}_{n\varphi}\to\infty$ as 
$n\to+\infty$. Given any large $n$, specify $k$ by the relation 
$n\in[N_k,N_{k+1})$. Let $D^k_j$ be any disk in $\scriptf(n)$. 
Consider the annular region $\scripta^k_j = \{z: 
\rho_k<|z-z^k_j|<\e_k\}$. This region is disjoint from $D^l_i$ for 
all $l\le k$, except for $(l,i)=(k,j)$. Thus in $\scripta^k_j$, 
$\Delta\varphi\equiv \sum_{l>k}\sum_i \Delta\varphi^l_i$. 
 
Define 
$\tilde\varphi_k = \sum_{l\le k}\sum_i \varphi^l_i$. 
We have $\int_{D^k_j}n\Delta\tilde\varphi_k 
= n\mu^k_j 
\in [\tfrac14,\tfrac12]$, and hence for any test function 
$u\in C^\infty_0(\Omega_0)$, 
\eqref{annularlowerbound} gives 
\[ 
\int_{\scripta^k_j} |u(z)|^2 \le 16\e_k^2 \int_{\scripta^k_j} 
|D_{n\tilde\varphi_k} u|^2. 
\] 
Let $E^k_j$ be the set of all radii $r\in[\rho_k,\tfrac12 \e_k]$ 
for which the circle $\Gamma_r=\{z: |z-z^k_j|=r\}$ intersects some 
closed disk $\ov{D^l_i}$ with $l>k$. The Lebesgue measure of 
$E^k_j$ is $\le 4\sum_{l>k} \e_l^{-2}\rho_l = 4\sum_{l>k} B^{2l} 
e^{-\sigma_l lB^{2l}}\ll B^{-k} = \e_k$, provided $k$ is sufficiently 
large. Moreover, for any $r\in[\rho_k,\e_k]$, 
\[ 
\int_{|z-z^k_j|<r} n\Delta (\varphi-\varphi^k_j) 
\le  n\sum_{l>k}\sum_i \mu^l_i. 
\] 
The sum over all $l\ge k+2$  contributes at most 
$4n\e_{k+2}^{-2}/N_{k+2} \le 4B^{2(k+2)}2^{-(B^{k+1}-B^k)} <1/8$ 
for all sufficiently large $k$, since $n<N_{k+1}$. By 
construction, any point $z^{k+1}_i$ satisfies 
$|z^{k+1}_i-z^k_j|\ge B^{-k-1} = B\rp \e_k$. Thus 
$\varphi^{k+1}_i$ contributes nothing to the integral, provided 
that $r\le \tfrac12 B\rp\e_k$. We therefore conclude that for any 
$r\in[\rho_k,\tfrac12 B\rp \e_k]$ and any $m\ge k$, 
\[ 
\tfrac14\le\int_{|z-z^k_j|<r} 
n\Delta\varphi^k_j \le \int_{|z-z^k_j|<r} n\Delta\tilde\varphi_m \le 
\int_{|z-z^k_j|<r}n\Delta\varphi^k_j + n\sum_{l\ge k+2}\sum_i \mu^l_i 
\le \tfrac12 + \tfrac18 
\le 
\tfrac58, 
\] 
provided as always that $k$ is sufficiently large. 
 
Therefore for any $r\in[\rho_k,\tfrac12 B\rp\e_k]\backslash E^k_j$, 
\begin{equation*} 
\int_0^{2\pi} |u(z^k_j+re^{i\theta})|^2 \,d\theta \le 16 r^2 
\int_0^{2\pi} |D_{n\tilde\varphi_m}u(z^k_j+re^{i\theta})|^2 
\,d\theta 
\end{equation*} 
by \eqref{circlelowerbound}.  Hence 
\[ 
\int_{r\in[\rho_k,\tfrac12 B\rp\e_k]\backslash E^k_j} 
\int_0^{2\pi} |u(z^k_j+re^{i\theta})|^2 \,d\theta\,rdr \le 4 
B^{-2}\e_k^2 \int_{\scripta^k_j} |D_{n\tilde\varphi_m}u(z)|^2. 
\] 
Since $\|\nabla\tilde\varphi_m-\nabla\varphi\|_{L^2} \to 0$, 
we may now conclude, 
by letting $m\to\infty$, 
that for each $n\in [N_k,N_{k+1})$, 
for each $j$ such that $D^k_j\in\scriptf(n)$, 
\begin{equation} \label{pivotalconclusion} 
\int_{r\in[\rho_k,\tfrac12 B\rp\e_k]\backslash E^k_j} 
\int_0^{2\pi} |u(z^k_j+re^{i\theta})|^2 \,d\theta\,rdr \le 
4B^{-2}\e_k^2 \int_{\scripta^k_j} |D_{n\varphi}u(z)|^2. 
\end{equation} 
 
We claim next that 
\begin{equation}\label{claim1} 
\int_{B(z^k_j, \e_k)}|u|^2 
\le 
C\e_k^2 \int_{B(z^k_j, \e_k)}\big|\nabla|u|\big|^2+ 
C \int_{r\in[\rho_k,\tfrac12 B\rp\e_k]\backslash E^k_j} 
\int_0^{2\pi} 
|u(z^k_j+re^{i\theta})|^2 
\,d\theta\,rdr, 
\end{equation} 
where $C$ is a constant depend only on $B$.  The proof of this 
claim follows from the Poincar\'e inequalities \eqref{poincare1} 
and \eqref{poincare2}, as follows. By \eqref{poincare1}, for any 
$r\in [\tfrac14 B^\rp \e_k, \tfrac12 B^\rp \e_k]$, 
\[ 
\int_{B(z^k_j, \tfrac14 B^\rp\e_k)}|u|^2 
\le 
C\e_k^2 \int_{B(z^k_j, \e_k)}\big|\nabla|u|\big|^2 
+C \e_k^2\int_0^{2\pi}|u(z^k_j+re^{i\theta})|^2 
\,d\theta. 
\] 
Integrating both sides with respect to $r$ over 
$ [\tfrac14 B^\rp \e_k, \tfrac12 B^\rp \e_k]\backslash E^k_j$, 
dividing both sides by $\e_k$, and using the fact that 
$|E^k_j|\ll \e_k$, we obtain 
\begin{equation}\label{claim2} 
\int_{B(z^k_j, \tfrac14 B^\rp\e_k)}|u|^2 
\le 
C\e_k^2 \int_{B(z^k_j, \e_k)}\big|\nabla|u|\big|^2+ 
C \int_{[\rho_k, \tfrac12 B^\rp \e_k]\backslash E^k_j} 
\int_0^{2\pi}|u(z^k_j+re^{i\theta})|^2\,d\theta\,rdr. 
\end{equation} 
Similarly, it follows from \eqref{poincare2} that 
for any $r\in [\tfrac18 B^\rp \e_k, \tfrac14 B^\rp \e_k]$, 
\[ 
\int_{\tfrac14 B^\rp\e_k<|z-z^k_j|<\e_k}|u|^2 
\le 
C\e_k^2 \int_{B(z^k_j, \e_k)}\big|\nabla|u|\big|^2 
+C \e_k^2\int_0^{2\pi}|u(z^k_j+re^{i\theta})|^2 
\,d\theta. 
\] 
Integrating both sides over $r\in 
[\tfrac18 B^\rp \e_k, \tfrac14 B^\rp \e_k]\backslash E^k_j$ 
and dividing by $\e_k$, we have 
\begin{equation}\label{claim3} 
\int_{\tfrac14 B^\rp\e_k<|z-z^k_j|<\e_k}|u|^2 
\le 
C\e_k^2 \int_{B(z^k_j, \e_k)}\big|\nabla|u|\big|^2+ 
C \int_{[\rho_k, \tfrac12 B^\rp\e_k] 
\backslash E^k_j}\int_0^{2\pi}|u(z^k_j+re^{i\theta})|^2 
\,d\theta\,rdr. 
\end{equation} 
Combining \eqref{claim2} and \eqref{claim3}, we then 
obtain \eqref{claim1}. 
 
Applying \eqref{poincare2} once more, we conclude that for 
any fixed finite constant $C'>0$, 
\begin{equation}\label{claim4} 
\int_{\e_k<|z-z^k_j|<C'\e^{1/2}_k} |u|^2 
\le C\e_k\log(1/\e_k) \int_{ \e_k<|z-z^k_j|<C'\e^{1/2}_k} \big|\nabla|u|\big|^2 
+ C\e_k \e_k^{-2}\int_{|z-z^k_j|\le\e_k} |u|^2 
\end{equation} 
By \eqref{claim1}, \eqref{claim4}, \eqref{pivotalconclusion}, and Lemma~\ref{Kato}, 
we have 
\begin{align*} 
\int_{B(z^k_j,C'\e_k^{1/2})} |u|^2 
&\le C\e_k\log(1/\e_k) \int_{B(z^k_j,C'\e_k^{1/2})} \big|\nabla|u|\big|^2 
+ C\e_k^{-1}\int_{B(z^k_j, \e_k)} |u|^2 
\\ 
&\le C\e_k\log(1/\e_k) \int_{B(z^k_j,C'\e_k^{1/2})} |D_{n\varphi}u|^2 
+ C \e_k \int_{B(z^k_j, \e_k)} |D_{n\varphi}u|^2 
\\ 
&\le C\e_k\log(1/\e_k) \int_{B(z^k_j,C'\e_k^{1/2})} |D_{n\varphi}u|^2. 
\end{align*} 
This holds for any $j$ such that $D^k_j\in\scriptf(n)$. 
 
It follows that 
\[ 
\int_{\cup_j B(z^k_j,C'\e_k^{1/2})} |u|^2 \le CB^{-k/2} 
\int_{\Omega_0} |D_{n\varphi}u|^2, 
\] 
where the union is taken over all $j$ such that $D^k_j\in\scriptf(n)$. 
Indeed, choose a maximal pairwise disjoint subfamily of all the balls 
$B(z^k_j,C'\e_k^{1/2})$, apply the preceding inequality to 
the tripled ball $B(z^k_j,3C'\e_k^{1/2})$ for each element of the subfamily, 
and sum. 
No point of $\complex$ belongs to more than a fixed number 
of the tripled balls, and each ball $B(z^k_j,C'\e_k^{1/2})$ 
is contained in at least one of the tripled balls. 
 
The construction guarantees that any point in the complement in 
$\Omega_0$ of $\cup_j B(z^k_j, 16\e_k^{1/2})$ either lies within 
distance $\e_k^{1/2}$ of the boundary of the unit ball $\Omega_0$, 
or is contained in $\cup_{1\le l<k}\cup_{i=1}^{m_l} D^l_i$. For 
the former region an application of the fundamental theorem of 
calculus, as in the proof of \eqref{poincare2}, gives the simple 
bound $\int_{1-\sqrt{\e_k}<|z|<1} |u|^2 \le 
C\sqrt{\e_k}\int|\nabla |u||^2$, which is dominated by 
$CB^{-k/2}\langle S_{n\varphi} u,\,u\rangle$ by Kato's inequality. 
$k=k(n)\to\infty$ as $n\to\infty$, and $\e_k\to 0$ as 
$k\to\infty$, so this region is satisfactorily under control. 
 
For any $n\in [N_k, N_{k+1})$ and any $l<k$,  let $\tilde\rho_l 
=(1-2/\log n)\rho_l$ and $\hat\rho_l =(1-3/\log n)\rho_l$.  Denote 
$\tilde{D}^l_i=B(z^l_i, \tilde\rho_l)$.  The electric potential 
term gives us 
\begin{equation} 
\sum_{l<k}\sum_i \int_{\tilde D^l_i} n\rho_l^{-2}\mu^l_i\cdot c_0 
n^{-1/2} |u|^2 \le \langle S_{n\varphi}u,\,u\rangle . 
\end{equation} 
Now 
\begin{equation*} 
n\rho_l^{-2} \mu^l_i 
 \ge n e^{2\sigma_l l B^{2l}}  N_{l+1}^{-1} 
= n e^{2\sigma_l l B^{2l}} 2^{-B^l} \ge n, 
\end{equation*} 
so we conclude that 
\[ 
\int_{\cup_{l<k}\cup_i \tilde D^l_i}|u|^2 \le Cn^{-1/2} \langle 
S_{n\varphi} u,\,u\rangle. 
\] 
Applying \eqref{poincare2} to the annulus $\{r<|z-z^l_i|<\rho_l\}$ 
for all $r\in [\hat\rho_l, \tilde\rho_l]$, then integrating (after 
multiplying both sides by $r$), we obtain 
\begin{align*} 
\int_{\tilde\rho_l<|z-z^l_i|<\rho_l} |u|^2 &\le 
(\rho_l^2-\hat{\rho}_l^2)\log(\rho_l/\hat\rho_l)\int_{D^l_i} 
\big|\nabla |u|\big|^2 
+\frac{2(\rho_l^2-\hat\rho_l^2)}{\tilde\rho_l^2-\hat\rho_l^2}\int_{\tilde 
D^l_i} |u|^2 \\ 
&\le C\big|\log(1-3/\log n)\big|\int_{D^l_i} \big|\nabla 
|u|\big|^2 + C\int_{\tilde D^l_i} |u|^2 . 
\end{align*} 
Once again $\big|\nabla|u|\big|$ may be replaced by 
$|D_{n\varphi}u|$ in the last integral, by Lemma~\ref{Kato}. Since 
every point of $\Omega_0$ belongs either to 
$\{1-\sqrt{\e_k}<|z|<1\}$, to some $D^l_i$ with $l<k$, or to 
$B(z^k_j, 16\e_k^{1/2})$ for some $j$ such that 
$D^k_j\in\scriptf(n)$, we conclude finally that 
\[ 
\int_{\Omega_0}|u|^2 \le C\min(B^{-k/2},n^{-1/2}, 
\big|\log(1-3/\log n)\big|) \langle 
S_{n\varphi} u,\,u\rangle, 
\] 
where $k=k(n)$ is determined by the relation 
$n\in[N_k,N_{k+1})$. 
 
It remains to prove that $\lim_{n\to\infty} \lambda^{\rm e}_{n\varphi} 
<\infty$.  This is easy.  Let $F_k$ be the functions constructed 
in the proof of Lemma~\ref{fine}.  Recall that $F_k$ is 
piecewise smooth, vanishing outside $\Omega_k$, 
$0\le F_k(z)\le 1$, $\|F_k\|_{L^2}\ge c>0$, 
and $\|\nabla F_k\|_{L^2}\le C<\infty$. 
Therefore 
\begin{equation*} 
\lambda^{\rm e}_{n\varphi} \le C\big( \int_{\Omega_k} |\nabla F_k|^2 + 
n\int_{\Omega_k} \Delta\varphi F_k^2\big)  
\le C \big(1 +n\sum_{l>k}\sum_j \mu^l_j\big)  
\le C\big(1+n\sum_{l>k} \nu_l\big). 
\end{equation*} 
This is $O(1)$ if $k$ is chosen $\ge n$, since $\sum_k k\nu_k<\infty$.
We thus deduce that $\lambda^{\rm e}_{n\phi}\le C<\infty$.  This concludes 
the proof of Theorem~\ref{main}. \qed 
 
\begin{remark} 
An interesting discussion of the lowest eigenvalue $\lambda^{\rm m}_\psi$ of $S_\psi$ 
on multiply connected domains with finitely many holes, in the special 
case where the winding number  corresponding to each hole is congruent to $\tfrac12$ modulo $\integers$,
appears in the work \cite{helfferetal} of Helffer et.\ al, who 
lift the problem to a twofold covering  
surface on which each winding number belongs to $\integers$.  
It may very well be possible to 
show in this way that the lowest eigenvalue for $S_\psi$ is large in 
certain situations, for instance when there are a large number of holes 
which are fairly densely distributed, if the winding numbers are half-integers. 
But our construction requires consideration of the case where the winding
number varies over a $\frac1N$--dense subset of an interval $[\delta,1-\delta]$ modulo $\integers$,
where $N\to\infty$.
We have followed a direct analytic path, 
based on Kato's inequality and the magnetic effect expressed through 
Lemma~\ref{magneticlemma}, whereby it is clear that 
the local effect produced by a single hole 
may be quantified in terms of the distance from a winding number to $\integers$, 
rather than its Diophantine character. 
\end{remark} 

\begin{remark}
By refining the estimates of this section slightly one can carry out the
construction so that $\Delta\varphi\in L(\log L)^\delta$ for any $\delta<1$.
It appears that one can get $\Delta\varphi\in L(\log L)(\log \log L)^{-C}$
for some finite $C$, but we have not verified this in detail.
\end{remark}

\section{Ground state energies in the smooth case} \label{section:smoothcase}

In this section we prove Theorem~\ref{thm:smoothcase}.
We are given that $\Delta\varphi$ is H\"older continuous of order $\alpha>0$,
and that $\lambda_{n\varphi}^{\rm e}$ remains bounded as $n\to+\infty$.

Before embarking on the proof, we pause to explain the underlying issues.
For large positive $n$,
$\int |\nabla u|^2 + n \int |u|^2\Delta\varphi$ will be large
relative to $\int|u|^2$, unless $u$ is supported mainly where 
$\Delta\varphi$ is nearly zero.
magnetic field into finitely many components for each $n$. 
Let $W$ be any connected component of the open set $\{\Delta\varphi>0\}$.
For any $n$, $\Delta\varphi\big|_W$ gives rise to a magnetic field 
whose strength (if $W$ is simply connected) on the complement of $W$ is governed by
the distance from $\int_W\Delta\varphi$ to $2\pi\integers$; if this distance
is nearly zero modulo $2\pi\integers$ then this part of the magnetic
field should not account for much of a discrepancy between 
$\lambda_{n\varphi}^{\rm e}$
and 
$\lambda_{n\varphi}^{\rm m}$.
Since all that is relevant is $n\int \Delta\varphi$ modulo $2\pi\integers$,
fields created by different components $W$ can interfere destructively with
one another.
Moreover, even if the field due to $W$ is strong, it is strong only near $W$; 
all that is required for 
$\lambda_{n\varphi}^{\rm m}$.
to be not much larger than
$\lambda_{n\varphi}^{\rm e}$
is for there to exist some suitably large subregion of $\Omega$ on which 
the net magnetic field is not very strong.
Thus it can be advantageous in the analysis to group components $W$ into 
clusters. Moreover, 
two or more components separated
by narrow necks will tend to act like a single larger component, since
a Brownian particle is unlikely to pass through a narrow neck without
straying into one of the components bounding it.
Thus breaking the support of $\Delta\varphi$ into
its topological components is inefficient.
Figure 2 illustrates some of these points.
\begin{figure}[htbp] \centering
 \ \psfig{figure=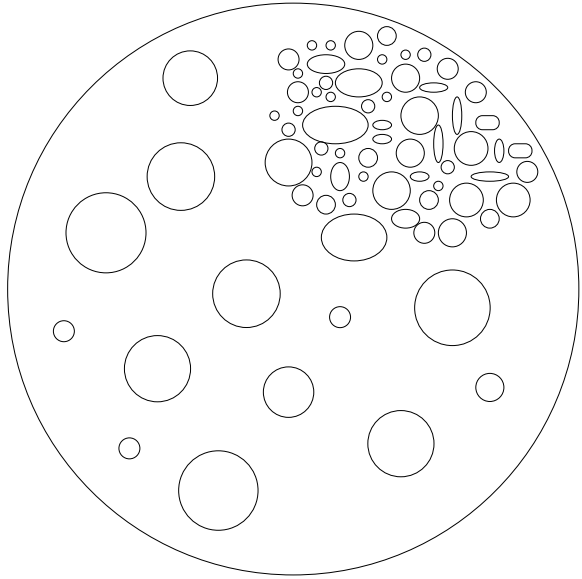,width=3in,height=3in}
\label{figure2}
 \caption{
A situation in which the support of the magnetic field
$\Delta\varphi$, represented by disks and ellipses, has many topological
components, which can be effectively organized into a small number of clusters.
}
\end{figure}

\begin{lemma} \label{trialfunctions}
Let  a subharmonic function $\varphi$ be given, with
$\Delta\varphi$ H\"older continuous of some order $\alpha>0$.
Suppose that $\lambda_{n\varphi}^{\rm e}$ remains bounded as $n\to+\infty$.
Then there exists $C<\infty$ such that
for any $\delta>0$ there exists a real-valued function $u\in C^\infty$
such that
$u$ is supported in $\Omega_0$,
$u\ge 0$, 
$\|\nabla u\|_{L^2}\le C$,
$\|u\|_{L^\infty}=1$, 
$\|u\|_{L^2}\ge C\rp$,
and $u(z)=0$ wherever $\Delta\varphi(z)\ge \delta$.
\end{lemma}

\begin{remark}
It follows from results in potential theory \cite{Fuglede99}, \cite{FuStraube02}
that there exists $u_0$ with $\nabla u_0\in L^2$,
supported in $\{z\in\Omega_0: \Delta\varphi(z)=0\}$,
so that $u_0\ne 0$ on a set of positive Lebesgue measure.
The desired function $u$ may be obtained by suitably mollifying $u_0$.
We have elected instead to give a self-contained proof.
\end{remark}

\begin{proof}
We are given that there exists $B<\infty$ such that for any $M<\infty$ there exists 
a $C^1$ function $v_M$ supported in $\Omega_0$ 
such that $\|\nabla v_M\|^2 + M\int |v_M|^2\Delta\varphi\le B^2$,
and $\|v_M\|_{L^2}\sim 1$.
By replacing $v_M$ by its absolute value we may assume that $v_M\ge 0$.
By replacing $v_M$ by $\max(v_M(z),c)-c$ for some sufficiently small $c>0$
we may make $v_M$ be supported in a compact subset of $\Omega_0$,
retaining uniform bounds on $v_M,\nabla v_M$.

Since $\|v\|_{L^4}$ is bounded by a constant times $\|\nabla v\|_{L^2}+ \|v\|_{L^2}$,
for large $\lambda>0$ we have $\int_{|v(z)|>\lambda} |v|^2
\le \lambda^{-2}\int_{|v(z)|>\lambda}|v|^4\le (B+1)C\lambda^{-2}$.
Therefore if we fix a sufficiently large constant $\lambda$ and  replace
$v_M$ by  $\min(v_M(z),\lambda)$, we still have 
$\|\nabla v_M\|^2 + M\int |v_M|^2\Delta\varphi\le B^2$
and $\|v_M\|_{L^2}\ge \tfrac12$, and have the additional property $\|v_M\|_{L^\infty}\le\lambda$
uniformly in $M$.

Choose a cutoff function $\eta=\eta_\delta\in C^\infty$, taking values in $[0,1]$, 
such that $\eta(z)\equiv 1$
wherever $\Delta\varphi(z)\le \delta/4$, and $\eta(z)\equiv 0$ wherever $\Delta\varphi(z)\ge\delta/2$.
Consider $u=u_{M,\eps}=\eta_\eps v_M$ where $M=M(\delta)$ is to be chosen sufficiently large.
$\|u\|_\infty$ is bounded above uniformly in $M$. $\|u\|_{L^2}$ is bounded below
by a strictly positive constant, provided that $M(\delta)$ is sufficiently large, 
since 
\[
\int_{\Delta\varphi(z)\ge\delta/4}|v_M(z)|^2
\le M\rp 4\delta\rp\int M\Delta\varphi|v_M|^2\le M\rp 2\delta\rp B^2
\]
may be made as small as desired by choosing $M$ sufficiently large.
Finally
$M\int |u|^2\Delta\varphi\le M\int |v_M|^2\Delta\varphi$ is uniformly bounded, while
\[
\|\nabla u\|_{L^2}\le \|\nabla v_M\|_{L^2} + \|v_M\nabla\eta\|_{L^2}
\le B + \|\nabla\eta\|_{L^\infty}\|v_M\|_{L^2(\text{support}\,(\nabla\eta))}.
\]
Since $\Delta\varphi\ge\delta/4$ on the support of $\nabla\eta$,
the last term is $O(M\delta\rp)^{-1/2}$, which tends to $0$ as $M\to\infty$.

The final step is to convolve with an approximation to the identity to produce a $C^\infty$
function; all the bounds continue to hold uniformly for a sufficiently fine approximation.
\end{proof}

Define $N_k = 2^k$ for each $k\in\naturals$; note that these differ from the
quantities denoted $N_k$ in previous sections.
For each $k\in\naturals$
let $u_k\in C^\infty$ satisfy the conclusions of Lemma~\ref{trialfunctions}
with $\delta = N_k^{-2}$. 
By Sard's theorem, there exist 
regular values of $u_k$ in $[\tfrac14,\tfrac12]$;
choose any such regular value and denote it by $c_k$.

Consider all connected components $V_j$ of $\{z\in\Omega_0: u_k(z)<c_k\}$.
Such a component is said to be harmless if
$\sup_{z\in V_j} \Delta\varphi(z)\le N_k^{-2}$, and
to be dangerous otherwise. 
Since $c_k$ is a regular value, there are only finitely many
such components, and each component of the boundary of any
$V_j$ is a smooth Jordan curve.


Let $\{W_i^k: 1\le i\le M_k\}$ be the collection of all
dangerous components $V_j$ of
$\{z\in\Omega_0: u_k(z)<c_k\}$.
$M_k$, the number of dangerous components, plays a central role in the analysis.


\begin{lemma}  \label{lemma:logcount}
For each $W_i^k$,
\[
\int_{W_i^k}|\nabla u_k|^2\gtrsim k\rp,
\]
uniformly in $i,k$.
Consequently
\[
M_k =  O(k).
\]
\end{lemma}

In contrast, if we were working with dangerous topological components
of $\{z: \Delta\varphi(z)>0\}$, the best bound would have roughly
the form $e^{ck}$. The size of $M_k$ will be the crucial element in
the pigeonhole argument of Lemma~\ref{lemma:pigeonhole}.

\begin{proof}
To verify this recall that $W_i^k$ contains a point $w_i$
for which $\Delta\varphi(w)\ge N_k^{-2}$, and that $W_i^k = V_j$ for some $j$.
The hypothesis of H\"older continuity implies that $\Delta\varphi(z)
\ge \tfrac12 N_k^{-2}$ for all $z$ in a disk $D_i$ of
radius $\gtrsim N_k^{-2/\alpha}$ centered at $w_i$.
Therefore $u_k(z)=0$ for all $z\in D_i$, and consequently
$D_i\subset V_j\subset W_i^k$.

Now consider the function $\tilde u$ which equals $u$ on $W_i^k$,
and equals $c_k$ on $\Omega^\dagger\setminus W_i^k$, where
$\Omega^\dagger$ denotes some fixed open ball which
contains the closure of $\Omega_0$. This function $\tilde u$
vanishes on the boundary of $D_i$, and equals $c_k\in[\tfrac14,\tfrac12]$
on the boundary of $\Omega^\dagger$.
Since $u_k\equiv c$ on the boundary of $W_i^k$,
$\nabla\tilde u\in L^2$ in the sense of distributions.
It follows from \eqref{poincare2} that
\[
\|\nabla\tilde u\|_{L^2}^2\gtrsim 1/\log(1/\rho)
\]
where $\rho$ denotes the radius of $D_i$.
Thus $\|\nabla\tilde u\|_{L^2}^2\gtrsim 1/\log(N_k)$.
But $\|\nabla u_k\|_{L^2(W_i^k)} = \|\nabla \tilde u\|_{L^2}$,
so the first conclusion is established.

The second conclusion follows directly. Since the sets $W_i^k$ are pairwise disjoint,
\[
M_k = \sum_{i=1}^{M_k}1
\lesssim  \sum_{i=1}^{M_k} k\int_{W_i^k}|\nabla u_k|^2
\le k\|\nabla u_k\|_{L^2}^2 =O(k).
\]
\end{proof}

The bound $M_k = O(k)$ is the best possible bound of this type, but is insufficient
for our purpose. The next lemma asserts an improvement for some subsequence.

\begin{lemma} \label{lemma:bettercount}
There exists a strictly increasing sequence $k_\nu\to\infty$
such that
\begin{equation}
M_{k_\nu}  \le \frac{k_\nu}{\log(k_\nu)\log\log(k_\nu)}.
\end{equation}
\end{lemma}

Our analysis requires a bound $M_{k_\nu}
= o(k_\nu/\log k_\nu)$; the factor of $\log \log k_\nu$
in the denominator serves to guarantee this but is not otherwise needed.

\begin{proof}
Let $K\in\naturals$ be large. Choose
$u_k=u$ and $c_k=c$ to be independent of $k$ for all $k\le K$;
these do still depend on $K$.
For each $2\le k\le K$ let $\{\tilde W_i^k: 1\le i\le m_k\}$
be the collection of all connected components of $\{z: u(z)\le c\}$
for which $\max_{z\in \tilde W_i^k}\Delta\varphi(z)\in [N_k^{-2}, N_{k-1}^{-2}) = [2^{-2k}, 2^{2-2k})$.
For $k=1$ the latter condition is instead $\max_{z\in \tilde W_i^k}\Delta\varphi(z)\ge \tfrac14$.
Then for any $k\le K$, $\{W_i^k\} = \cup_{l\le k} \{\tilde W^l_i\}$,
so $M_k = \sum_{l=1}^k m_l$.

Arguing as in the proof of Lemma~\ref{lemma:logcount} we find that
\[
\sum_{k=1}^K k\rp m_k \lesssim \sum_{k=1}^K\sum_{i=1}^{m_k} \int_{\tilde W_i^k}
|\nabla u|^2 \le \|\nabla u\|_{L^2}^2.
\]
uniformly in $K$.
Therefore $\sum_{k=1}^K k\rp m_k=O(1)$ uniformly in $K$.

By summation by parts, it follows that likewise
$\sum_{k=1}^K k^{-2}M_k=O(1)$, since the boundary terms $k\rp M_k$
remain uniformly bounded by Lemma~\ref{lemma:logcount}.
Therefore for any given $K'$, if $K$ is sufficiently large there exists $k\in[K',K]$ such that
$M_k\le k[\log k\cdot\log\log k]\rp$.
Applying this for a sufficiently rapidly increasing sequence of values of $K$ yields the lemma.
\end{proof}

Fix a sequence $(k_\nu)$ satisfying the conclusion
of Lemma~\ref{lemma:bettercount}. Henceforth we consider only indices
$k$ belonging to this sequence, but omit the subscript $\nu$
in order to simplify notation.
The possibly very sparse subsequence which yields the bounded limit infimum 
in Theorem~\ref{thm:smoothcase} is obtained via the following application
of the pigeonhole principle.

\begin{lemma} \label{lemma:pigeonhole}
Let $A\in[1,\infty)$ be sufficiently large.
Then for each sufficiently large $\nu$,
there exists $n_\nu\le N_{k_\nu}$ such that
\begin{equation}  \label{pigeonholebound}
\distance(n_\nu\int_{W_i^k}\Delta\varphi, 2\pi\integers)\le k_\nu^{-A}
\qquad\text{for all $1\le i\le M_k$.}
\end{equation}
\end{lemma}

\begin{proof}
Write $k=k_\nu$ and let $\eps_k = k^{-A}$.
Consider the torus $T_k=(\reals/2\pi\integers)^{M_k}$,
with one coordinate for each index $i\le M_k$.
In $T_k$ consider the sequence of points
$p_n$, where the $i$-th component of $p_n$ equals
$n\int_{W_i^k} \Delta\varphi$ modulo $2\pi\integers$.

Partition $T_k$ into $\lesssim \eps_k^{-CM_k}$
cubes of sidelength $\le\tfrac12\eps_k$.
$N_k/2$ is larger than the number of such cubes,
since
\[\eps_k^{-CM_k} \le k^{ACk/\log k\log\log k}
\le e^{ACk/\log\log k}\ll \tfrac12 N_k = 2^{k-1}\]
for all sufficiently large $k$
because of the extra factor of $\log\log k$.
Therefore by the pigeonhole principle, there exist
indices $1\le n'<n''\le N_k$ such that
$p_{n'},p_{n''}$ belong to the same cube.
Setting $n_\nu=n''-n'$, we conclude that $n_\nu$ has the desired property.
\end{proof}

\begin{lemma} \label{lemma:ntendstoinfty}
There exists a sequence of natural numbers $n_\nu\le N_{k_\nu}$
satisfying \eqref{pigeonholebound}, such that
$n_\nu\to\infty$ as $\nu\to\infty$.
\end{lemma}

\begin{proof}
Modify the proof of Lemma~\ref{lemma:pigeonhole} as follows: Let $(b_k)$ be
some nondecreasing sequence of natural numbers, 
such that $b_k\to\infty$ as $k\to\infty$, but 
$b_k/\log\log k\to 0$.
In the proof of Lemma~\ref{lemma:pigeonhole}, consider only
parameters $n$ which are integral multiples of $2^{b_k}$,
$b_k = b_{k_\nu}$.
The pigeonhole principle still applies,
provided that 
\[e^{ACk/\log\log k}< \tfrac12 N_k 2^{-b_k} = 2^{k-1-b_k}.\]
This holds for all sufficiently large $k$, since $b_k/ \log\log k
\to 0$.
Therefore there exists an index $n_\nu$ satisfying \eqref{pigeonholebound},
which is a positive integer multiple of $2^{b_{k_\nu}}$; in particular,
$n_\nu\ge 2^{b_k}$. Thus $n_\nu\ge 2^{b_{k_\nu}}\to\infty$.
\end{proof}

\begin{lemma} \label{lemma:gauge}
Let $\Omega_0\subset\complex$ be open, 
and let $U_1,\dots U_N$ be pairwise disjoint open subsets of $\Omega_0$,
all with smooth boundaries. 
Set $\Omega = \Omega_0\setminus\cup_j U_j$.
Suppose that $u$ is a real-valued harmonic function in $\Omega$
which has a multiple-valued real harmonic conjugate $v$ in $\Omega$
such that $e^{iv}$ is single-valued in $\Omega$.
Then for any function $\psi$ with $\nabla\psi\in L^2(\Omega)$ and $\Delta\psi\in L^1(\Omega)$,
the quadratic form $S_\psi$ is unitarily equivalent in $L^2(\Omega)$
to $S_{\psi-u}$. That is, there exists a unitary mapping $\scriptu$ on $L^2(\Omega)$
which preserves $C^1_0(\Omega)$
such that
$Q_{\psi-u}(f,f)= Q_\psi(\scriptu f,\scriptu f)$
for all $f\in C^1_0(\Omega)$.
\end{lemma}

\begin{proof}
$\scriptu f(z) = e^{iv(z)}f(z)$ does the job, as one sees via the relations
$v_y = u_x$, $v_x=-u_y$.
\end{proof}

\begin{lemma} \label{lemma:gauge2}
Let $\Omega_0,U_j$ be as in Lemma~\ref{lemma:gauge}, and suppose that they
are all simply connected.
Suppose that $u\in C^2(\Omega_0)$, that $u$ is harmonic in $\Omega = \Omega_0\setminus\cup_j U_j$,
and that $\int_{U_j}\Delta u\in 2\pi\integers$ for each index $j$.
Then $u$ has a multiple-valued real harmonic conjugate $v$ in $\Omega$
such that $e^{iv}$ is single-valued in $\Omega$.
\end{lemma}

\begin{proof}
Define $\tilde u_j$ to be the Newtonian potential of $\Delta\varphi\cdot\chi_{U_j}$.
Since $u-\sum_j\tilde u_j$ is harmonic in the simply connected domain $\Omega_0$,
it has a single-valued harmonic conjugate $v_0$ in $\Omega_0$.
Since $U_j$ is simply connected, 
the fundamental group of $\complex\setminus U_j$ is $\integers$, and hence
the condition $\int_{U_j}\Delta \tilde u_j=\int_{U_j}\Delta u\in 2\pi\integers$
guarantees that $\tilde u_j$ has a multiple-valued harmonic conjugate $v_j$
on $\complex\setminus U_j$ such that $e^{iv_j}$ is single-valued.
Hence $v=v_0+\sum_j v_j$ is a multiple-valued harmonic conjugate for $u$
in $\Omega_0\setminus\cup_j U_j$ such that $e^{iv}$ is single-valued.
\end{proof}

\begin{lemma} \label{lemma:aftergauge}
If $n_\nu$ satisfies \eqref{pigeonholebound} for each $\nu$, then
the magnetic ground state eigenvalues
$\sup_\nu{\lambda_{n_\nu\varphi}^{\rm m}}$ remain uniformly
bounded as $\nu\to\infty$.
\end{lemma}

\begin{proof}
Let $k=k_\nu$, and $n=n_\nu$.
To each set $W_i^k$ associate $W_i^{k*}$, the smallest open simply connected
set containing $W_i^k$; this equals the union of $W_i^k$ with all the
bounded connected components of $\complex\setminus W_i^k$.
It may happen that one $W_i^{k*}$ is properly contained in another; delete
all such $W_i^{k*}$ and retain only those which are maximal with respect to
inclusion. Reorder so that those which remain are denoted
$\{W_i^{k*}: i\le M_k^*\}$ where $M_k^*\le M_k$.

For each $1\le i\le M_k^*$ there is some disk $D_i\subset W_i^{k*}$
of radius $\rho_i \gtrsim N_k^{-2/\alpha}$. 
Let $\tilde D_i$ be the disk concentric with $D_i$, with half as large a radius.
Choose a function $h_i$ supported on
$\tilde D_i$ such that 
$\|h_i\|_{L^1}\lesssim k^{-A}$,
and $\int_{W_i^{k*}} (n\Delta\varphi-h_i)\in 2\pi\integers$.
Since $\distance(\int_{W_i^{k*}} n\Delta\varphi,\ 2\pi\integers)=O(k^{-A})$, 
such functions clearly exist.
Let $H_i$ be the Newtonian potential $(2\pi)\rp h_i*\log|z|$ of $h_i$.
Then because $\Omega_0$ is bounded and every point of $\Omega\setminus D_i$ lies
at a distance $\ge \tfrac12\rho_i$ from the support of $h_i$,
\[\|\nabla H_i\|_{L^2(\Omega_0\setminus D_i)}\le C(\log\rho_i)^{1/2}\|h_i\|_{L^1}
= O(k^{1/2} k^{-A}).\] 

By Lemmas~\ref{lemma:gauge} and \ref{lemma:gauge2},
$S_{n\varphi}$ is unitarily equivalent,
in $L^2(\Omega_0\setminus \cup_{i=1}^{M_k^*}W_i^{k*})$,
to $S_{\psi_k}$
where $\psi_k$ is the Newtonian potential of
$n\Delta\varphi\chi_{\Omega_0\setminus\cup_i W_i^{k*}} + \sum_i h_i$.

Consider the test function $f = u_k-c_k$ in $\Omega_0
\setminus\cup_i W_i^{k*}$, $f =0$ in $\cup_i W_i^{k*}$.
Since $u_k\equiv c_k$ on the boundary of each $W_i^{k*}$,
$\|\nabla f\|_{L^2}\le \|\nabla u_k\|_{L^2}$,
which is uniformly bounded. Moreover $\|f\|_{L^2}$ is bounded below
by a strictly positive constant, uniformly in $k$.

Since $f=O(1)$ in $L^\infty$,
it suffices to show that $S_{\psi_k}(f,f)$ is bounded above,
uniformly in $k$.
Since $\|f\|_{L^\infty}$ and $\|\nabla f\|_{L^2}$ are uniformly
bounded, it suffices to show that $\|\nabla\psi_k\|_{L^2}$
remains bounded as $k\to\infty$.

Now $n\Delta\varphi\chi_{\Omega_0\setminus\cup_i W_i^{k*}}=O(2^{-k})$
in $L^\infty$ norm, by the definition of harmless components $V_j$,
so the gradient of its Newtonian potential is $O(2^{-k})$ in $L^\infty$
and hence also in $L^2$.
We have already noted that
the gradient of the Newtonian potential $H_i$ of $h_i$ is
$O(k^{-A+1})$ in $L^2$ norm on the complement of $D_i$,
hence on the complement of $W_i^{k*}$.
There are $M_k^*\lesssim k$ indices $i$,
so in all, $\nabla\psi = O(k^{-A+2})$ in $L^2$ norm.
By choosing $A>2$ we can ensure that
this is $O(k\rp)$.
Thus $\|\nabla\psi_k\|_{L^2}\to 0$ as $k\to\infty$.
\end{proof}

Since $n_\nu\to\infty$ as $\nu\to\infty$, Lemma~\ref{lemma:aftergauge}
implies Theorem~\ref{thm:smoothcase}.

\begin{remark}
The above arguments actually prove that for a subharmonic
function $\varphi$ such that $\Delta\varphi$ is H\"older
continuous of some positive order, $\lim_{n\to\infty}
\lambda^{\rm e}_{n\varphi}=\liminf_{n\to\infty}\lambda^{\rm m}_{n\varphi}=\lambda$,
where $\lambda$ is the first eigenvalue of the Dirichlet Laplacian
of the fine interior of $\{\Delta\varphi=0\}$.  This is consistent
with the well-known phenomenon that a strong magnetic field creates
a Dirichlet boundary condition in the semi-classical limit.
\end{remark}

\begin{remark}
Theorem\ref{thm:smoothcase} also holds in the following slightly more general
form.  Let $A=(a_1, a_2)\in C^{1+\alpha}(\Omega_0, \R^2)$ and let
$V=|\partial a_1/\partial y -\partial a_2/\partial x|$.  Let
$H_A=-(\nabla -iA)^2$ and $H^0_V=-\Delta +V$.  Let
$\lambda^m_{tA}$ and $\lambda^e_{tV}$ denote respectively the
first eigenvalues of the Dirichlet realizations of $H_{tA}$ and
$H^0_{tV}$.  Then $\lim_{t\to+\infty} \lambda^m_{tA}=+\infty$ if and
only if $\lim_{t\to+\infty} \lambda^e_{tV}=+\infty$.  This can be
proved by observing that after a gauge transformation, one may
assume that $A=(-\varphi_y, \varphi_x)$ and $V=|\Delta \varphi|$
for some $\varphi\in C^{2+\alpha}(\Omega_0)$.  Details are left to
the interested reader.
\end{remark}

\section{Reduction to a problem in $\complex^1$}\label{reduction}
The goal of this section is to reduce questions of compactness and
property ($P$) on a complete Hartogs domain in $\C^2$ to questions
concerning semi-classical limits of Schr\"{o}dinger operators in
$\C^1$. Theorem~\ref{compactness} is then an easy consequence of
this reduction and Theorem~\ref{main}.

\begin{proposition}\label{reductionthm}
Let $\Omega=\{(z, w)\in\C^2; \ |w|<e^{-\psi(z)}, \ z\in\Omega_0\}$
be a complete Hartogs domain with smooth, strictly pseudoconvex
boundary near $b\Omega\cap\{w=0\}$.  Suppose that $\psi$ is a
continuous subharmonic function on $\Omega_0$ such that 
$\nabla\psi\in L^2_{loc}(\Omega_0)$ in the sense of distributions.
Then
\begin{enumerate}
\item If $b\Omega$ satisfies property ($P$) then
$\lambda^{\rm e}_{n\psi}(\Omega_0)\to \infty$ as $n\to\infty$.
If moreover $\Delta\psi$ is lower semicontinuous,
then the converse also holds.

\item The Kohn Laplacian has compact resolvent if and only
$\lambda^{\rm m}_{n\psi}(\Omega_0)\to \infty$ as $n\to\infty$.
\end{enumerate}
\end{proposition}

Note that $\Delta\psi$ is indeed lower semicontinuous
in the construction underlying Theorems~\ref{compactness}
and \ref{main}.

When $b\Omega$ is smooth, the necessity in
Theorem~\ref{reductionthm}~(2) was first established by
Matheos~\cite{Matheos97} and the other assertions of the lemma were proved
in~\cite{FuStraube02}. Although their results were stated only for
the $\dbar$-Neumann Laplacian, their proofs contains the proof for
the Kohn Laplacian as well.  Since here we have only minimal
assumption on regularity of the boundary, some
modifications are needed. We provide details for the reader's convenience.

We first prove the forward implication in (1).  Choose relatively compact
subdomains $\Omega_2$ and $\Omega_1$,
$\Omega_2\subset\subset\Omega_1$, such that $b\Omega$ is strictly
pseudoconvex over $\Omega_0\setminus\Omega_2$.

For any $M>0$, there exist a neighborhood $U$ of $b\Omega$ and
$g\in C^\infty(U)$ such that $-1\le g\le 0$, and $\partial\dbar
g\ge M$ on $U$.   Replacing $g$ by $\int_0^{2\pi} g(z,
e^{i\theta}w)d\theta$ if necessary, we may assume that $g(z, e^{i\theta}w)\equiv g(z,w)$
for all $z,w$ and all $\theta\in\reals$.
Let $\chi\in C^\infty_0(B(0,
1))$ be a Friedrich mollifier ({\it i.e.}, $\chi\ge 0$,
$\chi(z)=\chi(|z|)$, and $\int_\C \chi =1$).  Let $\chi_\delta
=(1/\delta^2)\chi(z/\delta)$ and
$\psi_\delta=\chi_\delta\ast\psi$.  Then $\{\psi_\delta\}$ is a
decreasing sequence of smooth subharmonic functions converging
locally uniformly to $\psi$ as $\delta\to 0^+$.   Furthermore,
$\psi_\delta\to\psi$ in $L^1_{loc}(\Omega_0)$.

For $z\in\Omega_1$ and sufficiently small $\delta$, let
$h_\delta(z)=g(z, e^{-\psi_\delta (z)})$.  A straightforward
calculation using polar coordinates $w=re^{i\theta}$ yields that
\[
(h_\delta(z))_{z\bar z}=\partial\dbar g(X, \ov{X}) -
(\psi_\delta)_{z\bar z}g_r e^{-\psi_\delta}
\]
where $X=(1, -2(\psi_\delta)_z e^{i\theta-\psi_\delta})$.  Let
$C_M=\max\{|g_r e^{-\psi}|; \ \ z\in\ov{\Omega}_1\}$; $C_M$ depends
on $M$ since $g$ does.  Then $-1\le
h_\delta \le 0$ and $\Delta h_\delta\ge 4M-C_M\Delta \psi_\delta$
on $\ov{\Omega}_1$.

Let $\eta\in C^\infty_0(\Omega_1)$, $0\le \eta\le 1$, and $\eta=1$
on $\ov{\Omega}_2$.  For any $f\in C^\infty_0(\Omega_0)$,
substituting $u=f\eta$, $b=h_\delta$, and $\psi=0$ into
(\ref{twistorinequality}), we have 
\[ 
\int_{\Omega_0} |(f\eta)_z|^2 \ge \frac{1}{e}\int_{\Omega_0}
(h_\delta)_{z\bar z} |f\eta|^2 .
\] 
Therefore, 
\begin{align*} 
\int_{\Omega_0} (|\nabla f|^2 + n\Delta\psi |f|^2) &\ge 
\int_{\Omega_0} (\frac1{2} |\nabla (f\eta)|^2 -|f\nabla\eta|^2 
+n\Delta\psi |f|^2) \\ 
&\ge \frac1{8e}\int_{\Omega_0} \Delta h_\delta  
|f\eta|^2 +\frac{n}{2}\int_{\Omega_0} \Delta\psi |f|^2 \\ 
&\ge \frac{M}{2e}\int_{\Omega_0} |f|^2 +\int_{\Omega_0}
(\frac{n}{4}\Delta\psi -\frac{C_M}{2e}\Delta\psi_\delta)|f|^2 
\end{align*} 
when $n$ is sufficiently large.  Letting $\delta\to 0$, we obtain
that 
\[
\int_{\Omega_0}(|\nabla f|^2 +n\Delta\psi |f|^2) \ge 
\frac{M}{2e}\int_{\Omega_0} |f|^2 +\int_{\Omega_0} 
(\frac{n}{4}-\frac{C_M}{2e}) \Delta\psi |f|^2 
\ge 
\frac{M}{2e}\int_{\Omega_0} |f|^2, 
\] 
provided that 
$n$ is sufficiently large relative to $M$.  Therefore
$\lambda^{\rm e}_{n\psi}(\Omega_0)\ge M/2e$ for all sufficiently
large $n$, or equivalently, 
$\lim_{n\to\infty} \lambda^{\rm e}_{n\psi}(\Omega_0) =\infty$. 
 
Now we prove the converse direction in (1).  It suffices to prove that
$K=b\Omega\cap\{(z, w)\in\C^2; \ z\in\ov{\Omega}_2\}$ satisfies
property (P).   Let $V=\{z\in\Omega_0; \ \Delta\psi>0\}$; $V$ is open
since $\Delta\psi$ is assumed to be lower semicontinuous. 
Let  $K_0=\overline{\Omega_0}\setminus V$.   Then $K_0$ is a compact subset of
$\overline{\Omega_0}$ and $\Delta\psi \equiv 0$ on $K_0$.   We claim
that $K_0$ has empty fine interior.  Otherwise, $K_0$ supports a
function $\xi \in W^1$ which is nontrivial, that is, is nonvanishing on some
set of positive Lebesgue measure (see \cite{Fuglede99}).  It
follows that $\lambda^{\rm e}_{n\psi}(\Omega_0)\le 
\|\nabla\xi\|^2/\|\xi\|^2<\infty$, which contradicts our 
assumption. By Proposition~1.11 in \cite{Sibony87}, $K_0$
satisfies property ($P$).  Therefore, for any $M>0$, there exist 
a neighborhood $U_0$ of $K_0$ and $b\in C^\infty(U_0)$ such that
$0\le b\le 1$ and $\Delta b\ge M$.   Since $\Delta\psi$ is lower
semicontinuous and $\overline{\Omega_0}\setminus U_0$ is compact, 
there exists $\e_0>0$ such that $\Delta\psi\ge \e_0$ 
on $\Omega_0\setminus U_0$.   Let $g_\delta(z, 
w) =M(|w|^2 e^{2\psi_\delta}-1) +b(z)$.  Then when $\delta$ is
sufficiently small, $|g_\delta|\lesssim 1$ and $\partial\dbar
g_\delta \gtrsim M$ on $K$.  Hence $K$ satisfies property ($P$).

The proof of the necessity in (2) is easy: it suffices to plug
$u(z)e^{in\theta}$ into the compactness estimate
(\ref{2Dcompactnessestimate}) and use the fact that
$\|u(z)e^{in\theta}\|^2_{-1}\lesssim \|u\|^2/n^2$.

We now prove the sufficiency in (2) by establishing the
compactness estimate (\ref{2Dcompactnessestimate}). By assumption,
there exists $N_\epsilon>0$ such that when $n>N_\epsilon$,
\[
\|v\|^2\le \epsilon \|L_{n\psi} v\|^2, \quad \text{for all}\ v\in
C^\infty_0(\Omega_0).
\]
Taking the complex conjugate, we deduce that when $n<-N_\epsilon$,
\[
\|v\|^2\le \epsilon \|\bar L_{n\psi} v\|^2, \quad \text{for all}\
v\in C^\infty_0(\Omega_0).
\]
Therefore, when $|n|>N_\epsilon$
\[
\|v\|^2\le \epsilon (\|L_{n\psi} v\|^2+\|\bar L_{n\psi} v\|^2),
\quad \text{for all}\ v\in C^\infty_0(\Omega_0).
\]

For $u\in C^\infty_0 (\Omega_0\times T)$, write
$u=\sum_{n=-\infty}^\infty u_n(z) e^{in\theta}$. Then
\begin{align*}
\|u\|^2 &=\sum_{n=-\infty}^\infty \|u_n(z)\|^2 \le \sum_{|n|\le
N_\epsilon} \|u_n(z)\|^2 + \epsilon \sum_{|n|>N_\epsilon}
(\|L_{n\psi} u_n \|^2 +\|\bar L_{n\psi} u_n\|^2 ) \\
        &=\epsilon(\|Lu\|^2+\|\bar L u\|^2) + \sum_{|n|\le N_\epsilon} \left(\|u_n\|^2
-\epsilon(\|L_{n\psi} u_n\|^2 + \|\bar L_{n\psi} u_n\|^2)\right).
\end{align*}
Since the inclusion $\{u\in L^2(\Omega_0); \ \|L_{n\varphi}u\|^2
+\|\overline{L}_{n\varphi} u\|^2\le 1\} \subset L^2(\Omega_0)$ is
compact(see \cite{AvronHerbstSimon78}), the last sum in the above
inequalities is less than or equal to $C_\epsilon\sum_{|n|\le
N_\epsilon} \|u_n\|^2_{-1}$ for some sufficiently large
$C_\epsilon$, depending only on $\epsilon$. The desired inequality
(\ref{2Dcompactnessestimate}) then follows from the fact that this
last sum is controlled by $\|u\|^2_{-1}$.

We now indicate the standard procedure to construct the
function $\psi$ in Theorem~\ref{compactness} from the function
$\varphi$ in Theorem~\ref{main}.   Let $\chi_1$ be a smooth 
function such that $0\le\chi_1\le 1$, $\chi_1=1$ for $t\le 4/3$,
and $\chi_1=0$ for $t\ge 3/2$.   Let $\chi_2$ be a smooth function
on $(0, \, 2)$ such that $\chi'_2$, $\chi''_2 >0$, $\chi_2=0$ for
$t\le 1$, and $\chi_2(t)=-\frac{1}{2}\log (4-t^2)$ when $t$ is
sufficiently closed to 2.  We may choose such $\chi_2$ so that its
second derivative is large enough on $[4/3, \, 3/2]$ to guarantee
that the Laplacian of $\psi(z)=\varphi(z)\chi_1(|z|)+\chi_2 (|z|)$
is strictly positive when $1<|z|<2$.
 
That $b\Omega$ does not satisfy property ($P$) is then a 
consequence of Proposition~\ref{reductionthm}~(1) and the fact
that $\lambda^{\rm e}_{n\psi}(B(0, 2))\le\lambda^{\rm e}_{n\varphi}(B(0,
1))\lesssim 1$.  It remains to prove that $\lambda^{\rm m}_{n\psi}(B(0, 
2))\to\infty$ as $n\to\infty$.   This is a consequence of the
combination of the facts that $\lambda^{\rm m}_{n\varphi}(B(0, 1))\to
\infty$, $b\Omega$ is strictly pseudoconvex over $1<|z|<2$, and 
$b\Omega\cap\{(z,\, w)\in\C^2; \ |z|=1\}$ satisfies property
($P$); the details of this argument follow below. 
 
Let $\chi_2$ be a smooth function supported in $(-1,\, 1)$ such
that $-1\le \chi_3\le 1$, $\chi_3=t$ for $-1/2\le t\le 1/2$.  For
any $M>0$, let $b(z)=\chi_3(M(|z|^2-1))-1$.   Applying
(\ref{twistorinequality}) with this choice of $b(z)$, we obtain
that for any $u\in C^\infty_0(B(0, 2))$, 
\begin{equation}\label{bigone}
\int_{B(0, 2)} |L_{n\psi} u|^2 \ge\frac{1}{e^2}\int_{B(0, 2)}
b_{z\bar z} |u|^2 \ge\frac{M}{e^2}\int_{||z|^2-1|<\frac1{2M}}
|u|^2 -C_M\int_{||z|^2-1|\ge\frac1{2M}} |u|^2 .  
\end{equation}
where $C_M$ is a constant depend only on $M$.  Let 
$\zeta(z)=(|z|^2-1)^2$.  Then 
\begin{align}\label{bigtwo} 
\int_{B(0, 1)} |L_{n\psi} u|^2 &\ge \int_{B(0, 1)} |\zeta  
L_{n\psi} u|^2 \ge \frac1{2}\int_{B(0, 1)} |L_{n\psi} (\zeta 
u)|^2-\int_{B(0, 1)} |u\frac{\partial\zeta}{\partial z}|^2 \notag\\
&\ge \frac1{2}\lambda^{\rm m}_{n\varphi}(B(0, 1))\int_{B(0, 1)} |\zeta 
u|^2 -4\int_{B(0, 1)} |u|^2 . 
\end{align} 
Note also that  
\begin{equation}\label{bigthree} 
\int_{B(0, 2)} |L_{n\psi} u|^2 \ge n\int_{1<|z|<2} \Delta\psi 
|u|^2 . 
\end{equation} 
Combining (\ref{bigone}), (\ref{bigtwo}), (\ref{bigthree}), and
the facts that $\lambda^{\rm m}_{n\varphi}(B(0, 1))\to\infty$, 
$\Delta\psi>0$ on $1<|z|<2$, and $\Delta\psi\to\infty$ as $|z|\to
2$, we obtain
\[
\int_{B(0, 2)} |L_{n\psi} u|^2 \ge \frac{M}{2e^2}\int_{B(0, 2)}
|u|^2
\]
when $n$ is sufficiently large.  Therefore $\lambda^{\rm m}_{n\psi}(B(0,
2))\to\infty$ as $n\to\infty$.

\section{Proof of Theorem~\ref{thm:dbarsmoothcase} } \label{section:lastreduction}

The equivalence of assertions (1) and (2) of Theorem~\ref{thm:dbarsmoothcase}
has been established by Matheos
\cite{Matheos97} and the implication (3) $\Rightarrow$ (1) is a consequence of Catlin's
theorem \cite{Catlin84}.  We need only show that (1) implies (3). This will be a 
consequence of Theorem~\ref{thm:smoothcase}; the proof of the reduction  is divided into two lemmas.
Let $\Omega=\{(z, w)\in\C^2;\  \rho(z, w)<0\}$ where $\rho$ is a
smooth defining function that is invariant under rotations in $w$.
Let
\begin{align*}
S_0&=\{(z, w)\in b\Omega; \ \ \frac{\partial\rho}{\partial w}(z,
w)=0\}; \\
S_k&=\{(z, w)\in b\Omega; \ \ |\frac{\partial\rho}{\partial w}(z,
w)|\ge 1/k\}.
\end{align*}
Then $b\Omega=S_0\cup(\cup_{k=1}^\infty S_k)$.  By Proposition 1.9
in \cite{Sibony87}, it suffices to prove that each $S_k$, $k=0, 1, \ldots$,
is $B$-regular.

\begin{lemma} 
If the $\bar\partial$-Neumann Laplacian $\square$ has compact resolvent
in $L^2(\Omega)$ then
$S_0$ is $B$-regular.\end{lemma}

\begin{proof}
Let $\pi\colon S_0\to \C$ be the projection to the $z$-plane. Let
$\widehat S_0=\pi(S_0)$.  According to Proposition 1.10 in \cite{Sibony87},
it suffices to prove that all fibers $\pi^{-1}(z_0)$, $z_0\in
\widehat S_0$, as well as $\widehat S_0$ itself, are $B$-regular.

We identify $\pi^{-1}(z_0)$ with its projection to the $w$-plane.
Note that $\pi^{-1}(z_0)$ is a union of circles centered at the
origin.  Since $b\Omega$ is variety-free, $\pi^{-1}(z_0)$ must
have empty fine interior, and hence is $B$-regular. Otherwise,
suppose $w_0$ is a fine interior point of $\pi^{-1}(z_0)$.  Then
$\pi^{-1}(z_0)$ contains a circle centered at $w_0$ (cf.\ \cite{Helms69},
Theorem 10.14), which implies that $b\Omega$ contains an annulus
(or a disc when $w_0=0$). But the presence of a complex variety
in $b\Omega$ forces $\square$ to have a noncompact resolvent
(cf.\ \cite{FuStraube01}), which is a contradiction.

It remains to prove that $\widehat S_0$ is $B$-regular.  In fact,
we will prove that $\widehat S_0$ has zero Lebesgue measure. Let
$(z_0, w_0)\in S_0$.  Since $\rho_z(z_0, w_0)\not=0$, 
we may assume without loss of generality that $\rho_y(z_0, w_0)\not=0$. Then
in a neighborhood $U$ of $(z_0, w_0)$, $b\Omega$ is given by
$y=\tilde\rho(x, w)$ where $\tilde\rho$ is rotation-invariant with respect to
$w$.   It follows that locally
\[
\pi(\ov{U}\cap S_0)=\{x+i\tilde\rho(x, |w|)\in\complex^1: \ \
\frac{\partial\tilde\rho}{\partial |w|}(x, |w|)=0\},
\]
which by Sard's theorem has zero Lebesgue measure.

\end{proof}

\begin{lemma} \label{lemma:2ndBregularity} Each $S_k$, $k=1, 2, \ldots$, is 
$B$-regular.
\end{lemma}

\begin{proof}It suffices to prove that for any $(z_0, w_0)\in S_k$, there
exists a neighborhood $U$ of $(z_0, w_0)$ such that $\ov{U}\cap
S_k$ is $B$-regular.

Since $|\rho_w(z_0, w_0)|=|\rho_{|w|}(z_0, |w_0|)|/2\ge 1/k$,
$b\Omega$ is defined near $(z_0, w_0)$ by a graph of the form
$|w|=e^{-\varphi(z)}$.  Assume that $\rho_{|w|}(z_0, |w_0|)<0$.
(The other case is treated similarly.)  Then there exist $a, b>0$
such that
\[
\Omega\cap U_{a, b} =\{(z, w); \ \ z\in B(z_0, a),\
e^{-\varphi}<|w|<|w_0|+b\} .
\]
where $U_{a, b}=B(z_0, a)\times\{|w_0|-b<|w|<|w_0|+b\}$. 
The pseudoconvexity of $\Omega$ implies that $\varphi$ is superharmonic
on $B(z_0,a)$.
Shrinking $a$ if necessary, we may also assume that $e^{-\varphi}\in
(|w_0|-b/2, \ |w_0|+b/2)$ for all $z\in B(z_0, a)$.

For any $\beta\in C^\infty_0(B(z_0, a))$ and any positive integer
$n$, consider the $(0, 1)$-form
\[ 
u_n=\begin{cases} \beta(z) w^{-n} d\bar z, \quad |w_0|-b<|w|; \\
                   0, \quad \text{otherwise}.
\end{cases}
\]
Then $u_n\in C^\infty (\ov{\Omega})$ and $\dbar u_n=0$.   Since
the $\dbar$-Neumann operator has compact resolvent, the canonical
solution operator $S$ is likewise compact (cf.\ \cite{FuStraube01}). 
It follows that for any $\e>0$,
there exists $C_\e>0$ such that $f_n(z, w)=S(u_n)$ satisfies
\[ 
\|f_n\|^2_\Omega\le \e \|u_n\|^2_\Omega + C_\e\|u_n\|^2_{-1, 
\Omega} .
\]
Since $\|u_n\|^2_{-1, \Omega}\lesssim (1/n^2)\|u_n\|^2_{\Omega}$,
there exists $N_\e>0$ such that
\begin{equation}\label{dbarestimate}
\|f_n\|^2_\Omega\le 2\e \|u_n\|^2_\Omega, \quad n\ge N_\e.
\end{equation}
Since $\partial f_n/\partial\bar w=0$ for $(z, w)\in \Omega$ with
$z\in B(z_0, a)$,  $f_n(z, w)$ is holomorphic in $w$ on
$e^{-\varphi}<|w|<|w_0|+b$ for any fixed $z\in B(z_0, a)$.
Furthermore, $f_n(z, w)=g_n(z) w^{-n}$ where $\partial
g_n/\partial \bar z =\beta(z)$ on $\Omega\cap U_{a, b}$.  Note
that the left hand side of \eqref{dbarestimate} is

\begin{align}\label{LHS}
\|f_n\|^2_\Omega &\ge \int_{B(z_0, a)} |g_n(z)|^2 dxdy
\int_{e^{-\varphi}}^{|w_0|+b} r^{-2n+1} dr \notag \\
&=\frac{1}{-2(n-1)}\int_{B(z_0, a)} |g_n(z)|^2
\left(e^{2(n-1)\varphi}-(|w_0|+b)^{-2(n-1)}\right)dxdy,
\end{align}
and the right hand side of \eqref{dbarestimate} is

\begin{align}\label{RHS}
2\e\|u_n\|^2 &\le 2\e\int_{B(z_0, a)} |\beta(z)|^2
dxdy\int_{e^{-\varphi}}^A r^{-2n+1} dr \notag \\
&\le \frac{2\e}{-2(n-1)} \int_{B(z_0, a)} |\beta(z)|^2
\left(A^{-2(n-1)}-e^{2(n-1)\varphi}\right).
\end{align}
(Here we assume that $\Omega\subset\{|w|<A/2\}$ for some constant
$A>0$.) Since
\[
e^{2(n-1)\varphi} -(|w_0|+b)^{-2(n-1)}=e^{2(n-1)\varphi}\left(
1-\left( (|w_0|+b)^{-1}e^{-\varphi}\right)^{2(n-1)}\right)
\]
and
\[
(|w_0|+b)^{-1} e^{-\varphi} \le \frac{|w_0|+b/2}{|w_0|+b}<1,
\]
it follows from \eqref{dbarestimate}, \eqref{LHS}, and \eqref{RHS}
that for any $\beta\in C^\infty_0(B(z_0, a))$, there exists
$g_n(z)$ such that $\partial g_n/\partial\bar z =\beta(z)$ and
\[
\int_{B(z_0, a)} |g_n(z)|^2 e^{2(n-1)\varphi}\le 3\e \int_{B(z_0,
a)} |\beta(z)|^2 e^{2(n-1)\varphi}
\]
when $n\ge N_\e$. A duality argument then yields
that
\[
\int_{B(z_0, a)} |\alpha(z)|^2 e^{-2(n-1)\varphi} \le 3\e
\int_{B(z_0, a)} |\alpha_z(z)|^2 e^{-2(n-1)\varphi}, \quad \forall
\alpha\in C^\infty_0(B(z_0, a)), n\ge N_\e.
\]
Substituting $\alpha = ue^{(n-1)\varphi}$ and then replacing $n-1$ by $n$,
this becomes
\[
\int_{B(z_0, a)} |u(z)|^2  
\le 3\e
\int_{B(z_0, a)} |(\partial_z + n\varphi_z)u|^2
\quad \forall
u\in C^\infty_0(B(z_0, a)), n\ge N_\e.
\]
Now $\varphi$ is superharmonic; $-\varphi$ is subharmonic.
$-\partial_z -n\varphi_z = L_{-n\varphi}$, so
this last inequality is equivalent to $\lambda^{\rm m}_{-n\varphi}(B(z_0, a))\to\infty$
as $n\to +\infty$, which by Theorem~\ref{thm:smoothcase} implies that
$\lambda^{\rm e}_{-n\varphi}(B(z_0, a))\to\infty$. Therefore, by (the
proof of) Proposition 9.1 (1), $b\Omega\cap \ov{U}_{a, b}$
satisfies property ($P$).
\end{proof}

\begin{remark}
In the above proof we assume only that $b\Omega$ is of class $C^{2+\alpha}$,
which is needed to invoke Theorem~\ref{thm:smoothcase}. 
In the $C^\infty$ case, Lemma~\ref{lemma:2ndBregularity} could be proved more quickly by combining
Theorem~\ref{thm:smoothcase} with the equivalence, established by Matheos \cite{Matheos97},
between compactness in the $\bar\partial$-Neumann problem and boundary compactness in the sense
of \eqref{2Dcompactnessestimate}.
\end{remark}

\begin{remark}
If $\{\Delta\varphi=0\}=W$ is constructed as in Section~\ref{setconstruction},
then the conclusion of Theorem~\ref{thm:smoothcase} remains valid
whenever $\Delta\varphi$ is lower semicontinuous and in $L^p$ for some
$p>1$. 
\end{remark}


\begin{thebibliography}{99}
\bibitem{Matheos97} 
P.~Matheos, \emph{A {H}artogs domain with no analytic discs in 
the boundary for which the $\dbar$-{N}eumann problem is not 
compact}, UCLA Ph.D.\ thesis,  1997.
\end{thebibliography}
\end{document}